\documentclass[nonblindrev]{informs3}

\DoubleSpacedXI % Made default 4/4/2014 at request
\usepackage{endnotes}
\let\footnote=\endnote

\usepackage{natbib}
 \bibpunct[, ]{(}{)}{,}{a}{}{,}%
 \def\bibfont{\small}%
 \def\BIBand{and}%

\TheoremsNumberedThrough     % Preferred (Theorem 1, Lemma 1, Theorem 2)
\ECRepeatTheorems

\EquationsNumberedThrough    % Default: (1), (2), ...
\usepackage{multirow}
\usepackage{algorithm}

\newcommand{\reff}[1]{(\ref{#1})}

\makeatletter
\newenvironment{breakablealgorithm}
  {% \begin{breakablealgorithm}
   \begin{center}
     \refstepcounter{algorithm}% New algorithm
     \hrule height.8pt depth0pt \kern2pt% \@fs@pre for \@fs@ruled
     \renewcommand{\caption}[2][\relax]{% Make a new \caption
       {\raggedright\textbf{\ALG@name~\thealgorithm} ##2\par}%
       \ifx\relax##1\relax % #1 is \relax
         \addcontentsline{loa}{algorithm}{\protect\numberline{\thealgorithm}##2}%
       \else % #1 is not \relax
         \addcontentsline{loa}{algorithm}{\protect\numberline{\thealgorithm}##1}%
       \fi
       \kern2pt\hrule\kern2pt
     }
  }{% \end{breakablealgorithm}
     \kern2pt\hrule\relax% \@fs@post for \@fs@ruled
   \end{center}
  }
\makeatother

\begin{document}
\RUNAUTHOR{Luo et al.}
\RUNTITLE{Optimal Portfolio Deleveraging Problem}
\TITLE{Effective Algorithms for Optimal Portfolio Deleveraging Problem with Cross Impact}
\ARTICLEAUTHORS{%
\AUTHOR{Hezhi Luo}
\AFF{Department of Mathematics,   Zhejiang Sci-Tech University, Hangzhou,  Zhejiang 310018, China. \EMAIL{hzluo@zstu.edu.cn}} %, \URL{}}
\AUTHOR{Yuanyuan Chen}
\AFF{Corresponding author. Department of Finance and Insurance, Nanjing University, Nanjing, Jiangsu 210093, China. \EMAIL{yychen@nju.edu.cn}}
\AUTHOR{Xianye Zhang}
\AFF{Department of Mathematics,  Zhejiang Sci-Tech University, Hangzhou,  Zhejiang 310018, China. \EMAIL{201920102008@mails.zstu.edu.cn}}
\AUTHOR{Duan Li \textbf{(\textsc{deceased})}}
\AFF{School of Data Science, City University of Hong Kong, Hong Kong. \EMAIL{dli226@cityu.edu.hk}}
\AUTHOR{Huixian Wu}
\AFF{Department of Mathematics,  Hangzhou Dianzi University, Hangzhou,  Zhejiang 310018, China. \EMAIL{hxwu@hdu.edu.cn}}
}
\date{\today}

\ABSTRACT{
We investigate the optimal portfolio deleveraging (OPD) problem with permanent and temporary price impacts, where the objective is to maximize equity while meeting a prescribed debt/equity requirement. %For the first time in the literature,
We  take the real situation with cross impact among different assets into consideration. The resulting problem is, however, a non-convex quadratic program with a quadratic constraint and a box constraint, which is known to be NP-hard.
In this paper,
we first develop a successive convex optimization (SCO) approach for solving the OPD problem and show that the SCO algorithm  converges to a KKT point of its transformed problem. Second, we propose an effective global algorithm for the OPD problem, which integrates
the SCO method, simple convex relaxation  and a branch-and-bound framework,  to identify a global optimal solution to the OPD problem within a pre-specified $\epsilon$-tolerance. We establish the global convergence of our algorithm and estimate its complexity. We also conduct numerical experiments to demonstrate the effectiveness of our proposed algorithms with both real data and randomly generated medium- and large-scale OPD  instances.}

%\HISTORY{}
\maketitle
\section{Introduction}
In the financial market, leverage has been widely used by investors to significantly increase their returns of an investment. Although high leverage may be beneficial in boom periods, it may cause serious cash flow problems in recessionary periods. During the recessionary periods, investors often need to unwind their portfolios to reduce leverage. As pointed out by \cite{d2009leverage} and \cite{cont2016fire}, a sudden deleveraging of large financial portfolios has been recognized as a destabilizing factor in recent financial crises. As a result, how to deleverage large portfolios is an important question not only to investors, but also to the whole financial market.

 To deleverage a large portfolio is in fact to liquidate a large portfolio to meet the leverage requirement.  When liquidating a large portfolio, the price impact of trading must be taken into account.
\cite{Madhavan00} suggests that the price impact of trading can be decomposed into the temporary and permanent price impacts.
While the temporary price impact measures the instantaneous price pressure that is resulted from trading, the permanent price impact measures the
change in the equilibrium price that depends on the cumulative trading amount, and is
independent of the trading rate. See \cite{Carlin07} for a detailed discussion on the price impact.

%\chen{Market impact has been incorporated in many portfolio liquidation models (See e.g. \cite{Bertsimas98, Almgren00, Almgren03, Schied09a, Gatheral11, tsoukalas2017dynamic} and the references therein). }

\cite{Brown10} first study the optimal deleveraging problem with both permanent and temporary price impacts of trading.
With the assumption that both permanent and temporary price impact matrices are positive diagonal matrices and the convexity assumption that the temporary price impact is greater than one half of the permanent price impact for each asset in the portfolio, \cite{Brown10} derive an analytical optimal trading strategy for the optimal deleveraging problem.
 Empirical studies \citep{Holthausen90,Sias01,Schied09b}, however, show that permanent price impact may dominate in block transactions and hence the convexity assumption does not always hold.

 Without this convexity assumption, \cite{Chen14} study the optimal deleveraging problem considered in \cite{Brown10}, which reduces to a separable and non-convex quadratic program with quadratic and box constraints. They propose an efficient Lagrangian method to solve this non-convex quadratic program  and present conditions under which their method provides an optimal solution. Furthermore, \cite{Chen15} extend the Lagrangian method to the optimal deleveraging problem with non-linear temporary price impact.  However, \cite{Chen14, Chen15} still make the assumption that both permanent and temporary price impact matrices are diagonal in their models.
In other words, no cross impact among different assets is considered in \cite{Brown10} or \cite{Chen14,Chen15}.

%It is worth mentioning that no cross impact among the assets is considered in \cite{Brown10,Chen15,Chen14}.
%That is, both permanent and temporary price impact matrices are diagonal.
%In other words, none of the existing works on the optimal deleveraging problem takes the cross impact among different assets into consideration.
However, the existence of cross impact among different assets has been both empirically documented and theoretically justified. \cite{fleming1998information} show that portfolio rebalancing trades from privately informed investors can lead to cross impacts, even between assets fundamentally uncorrelated. \cite{kyle2001contagion} argue that correlated liquidity shocks due to financial constraints can lead to cross liquidity effects. \cite{Andrade08} show that a trade imbalance in one stock affects not only its own price, but also the prices of other stocks with their equilibrium model and data from the Taiwan Stock Exchange.  \cite{pasquariello2013strategic} empirically suggest that the strategical trading activity of sophisticated speculators may induce the cross impact, and find that the impact matrix could be asymmetric with negative non-diagonal elements in their empirical analysis with data from NYSE and NASDAQ. \cite{Benzaquen17} state that the cross-correlations among different stocks are, to a large extent, mediated by trades with the data set of 275 US stocks. Recently, \cite{lioptimal} apply genetic algorithm to the optimal deleveraging problem with non-diagonal symmetric impact matrices, without further discussion on i) how the non-diagonal impact elements affect the trading strategies, and ii) the convergence and complexity of their algorithm.

In this paper, we investigate how to solve the optimal deleveraging problem with cross impact, by considering general permanent and temporary price impact matrices, which could be not only non-diagonal, but also asymmetric. The resulting problem is then a non-separable and non-convex program with a quadratic objective function and non-homogeneous quadratic and box constraints, which is known to be NP-hard.
%In particular, our problem is a non-convex quadratic program with one quadratic and box constraints.

There have been quite a few works on the algorithm to solve a quadratic program with a quadratic constraint (QCQP). They, unfortunately, cannot be directly applied to solve our problem in this paper.
For example,  based on the celebrated $S$-lemma (cf. \cite{Polik07}), such a QCQP problem can be reformulated as a semi-definite programming problem (cf. \cite{Ben-Tal95,More93}).  \cite{Ben-Tal14} show that it can also be transformed into an equivalent second-order cone programming problem
when the quadratic forms are simultaneously diagonalizable. By exploiting the block separability of the canonical form, \cite{Jiang2018} show that all generalized trust
region subproblems with an optimal value bounded from below are second order cone programming (SOCP) representable.
% However, these methods do not deal with the QCQP with box constraints, and thus are not directly applicable to our problem.
However, these methods are not directly applicable to our problem, due to the box constraints.
\cite{Nesterov98} and \cite{Ye99} study the semi-definite relaxation method for certain quadratic programs with quadratic and/or box constraints.
 However, because  the quadratic constraint of our problem contains a linear term,  their method is also not directly applicable.  \cite{Audet00} propose a branch-and-cut algorithm for generic nonconvex QCQPs. \cite{Linderoth05} propose a simplicial branch-and-bound (B\&B) algorithm by reformulating the linear outer approximation of the non-convex constraints.
However, as reported in \cite{Audet00} and \cite{Linderoth05}, these two algorithms are targeted to solve small-scale instances. \cite{Burer08,Burer09} and \cite{Burer2012} propose the B\&B algorithms for non-convex quadratic programs with linear and/or box constraints by using  semi-definite relaxations of the first-order Karush-Kuhn-Tucker (KKT) conditions of the problem with finite KKT-branching.
More recently, \cite{lblp19} develop a new global algorithm for a non-convex quadratic program with linear and convex quadratic constraints by combining
 several simple optimization techniques  such as  the alternative direction method, the B\&B framework and the convex relaxation.
  %\cite{ldpjl20} further extended this approach to the worst-case linear optimization with uncertainties that arises in estimating the systemic risk in financial systems.
However, these algorithms cannot be directly applied  to solve our problem, due to the non-convex quadratic constraint in our problem.

Inspired by \cite{lblp19}, we develop an effective global algorithm for the optimal portfolio deleveraging (OPD) problem with cross impact in this paper. In particular, we first reformulate our problem as an equivalent D.C. (difference of convex functions) program by making use of the  spectral decomposition  and  the simultaneous diagonalization method of two semi-definite matrices (cf. \cite{Newcomb61}). We then propose a convex quadratic approximation for the transformed problem  via linearizing the concave quadratic terms in the objective function and the constraint of the transformed problem. Based on the solution to this convex approximation problem, we develop a successive convex optimization (SCO) approach for solving this transformed problem.
 We show that the solution sequence generated by our SCO approach converges to a KKT point of the transformed problem. By replacing the concave quadratic term in the transformed problem with its convex envelope, we also derive a simple convex relaxation for the transformed problem and estimate the gap between the relaxation and the transformed problem. Finally, we combine the SCO algorithm with the B\&B framework and the convex relaxation, to develop an efficient global algorithm (called SCOBB) that can find a globally optimal solution to the transformed problem within a pre-specified $\epsilon$-tolerance. We establish the convergence of the algorithm, and show that, the SCOBB algorithm has a complexity bound $ \mathcal{O} \left(N\prod\limits_{i=1}^{r}\left\lceil{\frac{\sqrt{r+\rho_1s}(z_u^i-z_l^i)}{2\sqrt{\epsilon}}}\right\rceil\right)$,
 where $N$ is the complexity to  solve the relaxed subproblem (a convex QP).

The contribution of this paper is thus twofold. Firstly, from the perspective
of the optimal deleveraging problem of a portfolio, we relax the non-realistic assumption of diagonal price impact matrices, which is usually made for mathematical and computational convenience, and derive analytical properties and numerical strategies on how to deleverage a portfolio with the consideration of cross impacts among different assets. With the non-diagonal cross impact matrices, we show that the asset, which is more liquid and less positively correlated with the other assets, would be prioritized for selling. Secondly, from the perspective of the solution to a non-convex quadratic program, we propose two efficient algorithms: the SCO algorithm and the SCOBB algorithm. We show that the SCO algorithm converges to a KKT point of the transformed problem,
%and thus serve the purpose as a good upper bound in the B\&B algorithm.
while the SCOBB algorithm  can find its global $\epsilon$-optimal solution  in $ \mathcal{O} \left(N\prod\limits_{i=1}^{r}\left\lceil{\frac{\sqrt{r+\rho_1s}(z_u^i-z_l^i)}{2\sqrt{\epsilon}}}\right\rceil\right)$ time.

Furthermore, we demonstrate the application of our algorithms to the optimal portfolio deleveraging problem with the historical data from NASDAQ in our numerical experiments. To illustrate the efficiency of our algorithms, we also randomly generate medium- and large-scale instances of the OPD problem to compare our two algorithms with the global optimization package BARON (cf. \cite{int:Sahinidis96}). We show that our SCOBB algorithm provides the globally optimal deleverage strategies with acceptable computational time for medium- and large-scale instances. Although we cannot prove that our SCO algorithm provides a global optimal solution, it can solve large-scale instances within several minutes, and always provides the global optimal solution in our numerical experiments.

Our work is also related to the literature on the portfolio execution problem.  \cite{Bertsimas99} develop an approximation algorithm to solve the multi-asset portfolio problem. \cite{Almgren00} briefly discuss the multi-asset portfolio problem with a mean-variance objective in their appendix and obtain a solution for the special case without cross impact. \cite{Tsoukalas17} analyze the optimal execution problem of a portfolio manager trading multiple assets with the cross impact.  The main difference between these works and ours is the motivation of the problem. While portfolio execution problem is to liquidate the whole portfolio, our optimal deleveraging problem is to liquidate some of the assets inside a portfolio to meet the leverage requirement.

The remaining of this paper is organized as follows. In Section~\ref{sec:ODP}, we  formulate the model of the OPD problem with the cross impact and derive analytical results on the optimal trading strategy. In Section~\ref{sec:formulation-SD}, we propose an equivalent reformulation for the OPD problem based on simultaneous diagonalization (SD) method. In Section~\ref{sec:SCO}, we propose the SCO algorithm for the transformed problem  and investigate its convergence properties. In Section~\ref{sect:GOA}, by integrating the SCO method with the SD-based convex relaxation, we develop a B\&B algorithm to find the global optimal solution to the transformed problem.  We establish the global convergence of the proposed algorithm and estimate its complexity.
  We test the performance of the proposed algorithms with historical data from NASDAQ and randomly generated instances in  Section~\ref{sec:Experiment}.  Finally, we present conclusions and future research directions in Section~\ref{sec:Conclusion}. The proofs of the technical results are all given in the e-companion to this paper.

\section{Optimal Portfolio Deleveraging Problem}\label{sec:ODP}
In this section, we  first present the model formulation for the optimal portfolio deleveraging problem, as in \cite{Brown10} and \cite{Chen14}. The main difference between our model and theirs is that we consider non-diagonal price impact matrices. We then explore some properties of the optimal trading strategy for our model.

\subsection{Problem Formulation}
  We consider a risk-neutral investor who trades a portfolio of $m$ assets in continuous time over a finite horizon. His initial position on this portfolio is $x_0\in{\mathbb R}^m$. The execution prices of the assets are modeled by
 \begin{equation}\label{eqn:price}
 p_t=q+\Gamma x_t +\Lambda y_t,
 \end{equation}
which is a multi-dimensional version of the pricing equation used in \cite{Carlin07}. Here, $p_t, x_t, y_t \in{\mathbb R}^m$ are the vectors of prices, holdings and trading rates of the $m$ assets at time $t$, respectively. The matrix $\Gamma=(\gamma_{ij})_{m\times m}$ documents permanent price impact coefficients of the cumulative trading amount $x_t$, where the element $\gamma_{ij}$ represents the permanent price impact parameter of asset $i$ on asset $j$. Similarly, the matrix $\Lambda=(\lambda_{ij})_{m\times m}$ is the temporary price impact matrix of trading associated with the trading rate $y_t$, whose element $\lambda_{ij}$ represents the temporary price impact parameter of asset $i$ on asset $j$. Both the initial holdings and initial prices are positive: $x_{0} >0$ and $p_{0}>0$. It is obvious that the trading amount and trading rate satisfy
$$x_t-x_0=\int_0^ty_sds.$$
In the same spirit as \cite{Brown10} and \cite{Chen14}, we assume a finite trading horizon of length $T=1$ and a trading strategy with a constant trading rate $y_t=y/T =y$, $t\in [0,1]$, where $y=x_1-x_0\in {\mathbb R}^m$ is the vector of cumulative trading amounts during the trading period.
Before trading, the asset prices are given by $p_0=q+\Gamma x_0$. After trading, the asset prices become $p_{1+}=q+\Gamma x_1=p_0 +\Gamma y$.
Because of the margin requirements imposed by lenders or the regulatory constraints, as pointed out by \cite{Brown10}, the trader has to guarantee the ratio of debt to equity to be no higher than a predetermined bound $\rho_1$. The equity is the difference between the value of the portfolio and the liability. In the beginning of the time horizon, the trader has an initial liability $l_{0}>0$ and an initial equity
$$e_{0} = p_0^Tx_0-l_0 > 0,$$
with the ratio of debt to equity exceeding $\rho_1$. In other words, we assume that the leverage requirement is not satisfied before trading, i.e., $\rho_1e_0<l_0$.
As a result, the trader has to liquidate some of his assets to meet the requirement. By trading a cumulative trading amount $y$ during the time horizon [0,1], the trader would generate the amount of cash as
\[\kappa(y)=\int_0^1-p_t^Ty_tdt=\int_0^1-\left(p_0+\Gamma\int_0^ty_sds+\Lambda y_t\right)^Ty_tdt=-p_0^Ty-y^T\left(\Lambda+\frac{1}{2}\Gamma\right)y.\]
The liability after trading can be obtained by subtracting
the cash amount expressed above from $l_0$:
\begin{equation}\label{eqn:liability}
l_1(y)=l_0-\kappa(y)=l_0+p_0^Ty+y^T\left(\Lambda+\frac{1}{2}\Gamma\right)y.
\end{equation}
%And the assets prices become $p_{1+}=p_0 +\Gamma y$.
Note that the asset prices after trading are $p_{1+}=p_0 +\Gamma y$. Thus the equity after trading is
\begin{equation}\label{eqn:equity}
e_1(y)=p_{1+}^Tx_1-l_1(y)=-y^T\left(\Lambda-\frac{1}{2}\Gamma\right)y+x_0^T\Gamma y+p_0^Tx_0-l_0.
\end{equation}
The trader  has to guarantee the leverage ratio after trading not exceeding  $\rho_1$, i.e.,
\begin{eqnarray}\label{eqn:leverage-constraint2}
l_1(y)/e_1(y)\le \rho_1~~ \Leftrightarrow~~ \rho_1e_1(y)-l_1(y)\ge 0.
\end{eqnarray}
Furthermore, with the same spirit as \cite{Brown10}, we are interested in modeling liquidity shocks that force the trader to quickly sell asset. As a result, the trader is restricted from increasing positions and from short selling during this fire sale. This corresponds to the box constraint $-x_0\le y\le 0$.
The problem is thus to find an optimal trading strategy that maximizes the equity $e_1(y)$ subject to both the leverage constraint and the box constraint. We formulate it as the following quadratic programming problem:
\begin{equation}\label{OPDP}
\begin{split}
\max\limits_{y\in{\mathbb R}^m}~& e_1(y)=-y^T\left(\Lambda-\frac{1}{2}\Gamma\right)y+x_0^T\Gamma y+e_0 \\
s.~t.~&  -y^T\left(\Lambda+\frac{1}{2}\Gamma+\rho_1\left(\Lambda-\frac{1}{2}\Gamma\right)\right)y-(p_0-\rho_1\Gamma x_0)^Ty-l_0+\rho_1e_0\ge0,\\
& -x_0\le y\le0.
\end{split}\tag{$P$}
\end{equation}
The first constraint corresponds to the leverage requirement \eqref{eqn:leverage-constraint2} with the liability and equity derived in \eqref{eqn:liability}-\eqref{eqn:equity}.

Because the prices of any two assets may interact with each other in real market, as we discussed in Section 1, the two matrices $\Gamma$  and $\Lambda$ could be non-diagonal and asymmetric. As a result, our model do not have any assumption on the price impact matrices. In addition, we assume that the leverage requirement can be always satisfied by liquidating all the assets:
\begin{assumption}\label{asmp1}
$l_1(-x_0)<0$. % and $\rho_1e_0-l_0<0$.
\end{assumption}

Note that the cumulative trading amount $y=-x_0$ implies that the investor liquidates all the assets. As a result, this assumption states that the investor should have no liability if he liquidates all his assets. With this assumption, our model focuses on the traders who do not need to liquidate all his assets to meet the liability obligation. If an investor cannot meet his liability obligation by liquidating all his assets, he would be bankrupted, and should just liquidate all his assets, instead of considering any optimization problem. Under this assumption, we show in the following proposition that if the assets are sold out,
then the investor's equity is positive, and that the leverage ratio is less than the predetermined level $\rho_1$.
\begin{proposition}\label{prop:property-problem}
$e_1(-x_0)>0$ and $\rho_1e_1(-x_0)-l_1(-x_0)>0$.
\end{proposition}
This proposition also implies that under Assumption~\ref{asmp1}, Problem (\ref{OPDP}) has a non-empty feasible strategy set, because the strategy $y=-x_0$ is feasible.

Note that we do not have any assumption on the price impact matrices $\Gamma$ and $\Lambda$. As a result, the problems with positive diagonal matrices in \cite{Brown10} and \cite{Chen14} are actually special cases of Problem (\ref{OPDP}). To deal with the difficulties introduced by the cross terms in the impact matrices, we present an equivalent reformulation for Problem (\ref{OPDP}) based on the spectral decomposition and simultaneous diagonalization method in Section~\ref{sec:formulation-SD}, and develop  algorithms in the subsequent sections. Before that, we first discuss the properties of our optimal strategy.

\subsection{Properties of Optimal Trading Strategy}\label{sec:property}
We first show that the leverage constraint of Problem (\ref{OPDP})
at the optimal solution is active under certain condition.
\begin{proposition}\label{prop-optimality}
If the price impact matrices $\Gamma$  and $\Lambda$ are non-negative and the permanent price impact matrix $\Gamma$ is symmetric,
then the leverage constraint of Problem (\ref{OPDP}) is active at its optimal solution.
\end{proposition}

 Proposition~\ref{prop-optimality} states that the optimal trading strategy precisely achieves the maximal allowed leverage ratio when the price impact matrices are non-negative and the permanent price impact matrix is symmetric. Under such a condition, it is suboptimal to further reduce the leverage
ratio because of the trading cost caused by market impact. Intuitively, during a fire sale, traders may hold the view that they should always, but not only under certain condition, liquidate the assets to precisely satisfy the leverage requirement, instead of selling more. However, if the price impact of one asset on another asset is negative, which is to say, there is hedging property inside the portfolio, selling one asset can improve the value of the other asset. The trader thus might be willing to sell more to improve his final equity. Note that, as we do not impose any assumption on the price impact matrices, the trader is  possibly willing to sell more in our model if some price impact parameters are negative.

 Proposition~\ref{prop-optimality} also implies that the positive diagonal impact matrices discussed in \cite{Chen14} and \cite{Brown10} can guarantee the activation of the leverage constraint.
 In addition, note that we only provide a sufficient, but not necessary condition for the activation of the leverage constraint.
 A symmetric permanent price impact matrix is thus not necessary for the activation of the leverage constraint. %In fact, we will show in our numerical experiments \chen{with real data} in Section~\ref{sec:Experiment} that the leverage constraint could be active for the example with asymmetric price impact matrices with several negative elements.
We show in Example~\ref{exam2} later that the leverage constraint could be active for Problem~(\ref{OPDP}) with asymmetric price impact matrices with several negative elements. This can be also verified by Example~\ref{exp_nasdaq} in our numerical experiments with real data in Subsection~\ref{sect:tprd}.

Next, we present a result on the asset trade priority, which are also discussed in \cite{Chen14} and \cite{Brown10}. We discuss this property without the assumption of diagonal price impact matrices. Define the average cross-stock impact matrices as
\begin{equation}\label{eqn:sys-matrices}
\hat\Gamma=(\hat\gamma_{lk})_{m\times m}:=\frac{1}{2}(\Gamma+\Gamma^T), \quad \hat\Lambda=(\hat\lambda_{lk})_{m \times m} :=\frac{1}{2}(\Lambda+\Lambda^T).
\end{equation}
The average cross-stock impact matrices describe the long-term and short-term co-movement between different assets' prices and trades. In particular, if $\hat\gamma_{lk}> \hat\gamma_{jk}$, then we can state that the mutual impact between assets $l$ and $k$ is higher than that between assets $j$ and $k$. As pointed out by \cite{Andrade08}, the mutual price impact is higher among assets with more positively correlated cash flows than among stocks with less positively correlated cash flows. We can thus understand these average cross-stock impacts as the indicator of the correlations for the cash flow of  the short-term or long-term assets. %Denote $p_{0}=(p_{0,1},\ldots,p_{0,m})^T$.
% the following proposition shows that trading priority is given to assets with smaller price impact.

\begin{proposition}\label{prop-priority}
Suppose assets $i$ and $j$ have the same initial
price and holding: $p_{0,i}=p_{0,j}$ and $x_{0,i}=x_{0,j}$. If the following four conditions hold, then the $i$th asset is prioritized for selling,
i.e., $y_i^*\le y^*_j$.
\begin{description}
\item[\rm (i)] %the price of the $i$th asset is less significantly affected by the trades of these two assets than the $j$th asset; meanwhile, the trades of the $i$th asset also less significantly affects the prices of these two assets than the $j$th asset, i.e.,
$\hat{\lambda}_{ii}\le \hat{\lambda}_{ij}=\hat{\lambda}_{ji} <\hat{\lambda}_{jj}$ and $\hat{\gamma}_{ii}\le \hat{\gamma}_{ij}=\hat{\gamma}_{ji} <\hat{\gamma}_{jj}$.
\item[\rm(ii)]For $\forall k\neq i,j$, %the average permanent cross-impact between $i$th asset and the $k$th asset is no higher than that between the $j$th asset and the $k$th asset, i.e.,
$\hat{\lambda}_{ik}\le \hat{\lambda}_{jk}$ and $\hat{\gamma}_{ik}\le \hat{\gamma}_{jk}$.
\item[\rm (iii)] $\gamma_{ii}-\gamma_{ij}\le\gamma_{jj}-\gamma_{ji}$.
\item[\rm(iv)] For $\forall k\neq i,j$, $\gamma_{ki}-\gamma_{ik} \le \gamma_{kj}-\gamma_{jk}$.
\end{description}

\end{proposition}

The condition (ii) in Proposition~\ref{prop-priority} implies that the $i$th asset is less positively correlated with the other assets. Meanwhile, the price of a more liquid asset is supposed to be less significantly affected by the liquidity shocks. The condition (i) in Proposition~\ref{prop-priority} thus implies that the $i$th asset is more liquid than the $j$th asset. The difference $\gamma_{jj}-\gamma_{ji}$ depicts to what extent the price of asset $j$ is more significantly affected by the trades of asset $j$ than the price of the $i$th asset. Similarly, the difference $\gamma_{kj}-\gamma_{jk}$ depicts to what extent the price impact of the trades of asset $k$ on asset $j$ is more significant than the price impact of the $j$th asset's trades on the asset $k$. As a result, the conditions (i)-(iv) in Proposition~\ref{prop-priority} imply that the $i$th asset is less positively correlated with the other assets and more liquid than the $j$th asset. Proposition~\ref{prop-priority} states that, if the $i$th asset is regarded as more liquid and less positively correlated with the other assets, than the $j$th
asset, it is prioritized for selling.

To illustrate the sufficiency and applications of the conditions in Propositions~\ref{prop-optimality}-\ref{prop-priority}, we construct the following two examples.

\begin{example}\label{exam1}{\rm ~~Consider an instance of Problem~(\ref{OPDP}) with $m=3$ and the following parameters}:
\begin{eqnarray*}
&&\Lambda=\left(\begin{array}{ccc} 0.0052&0.0083&0.0087\\0.0037&0.0085&0.0059\\
    0.0037&0.0094&0.0093\end{array}\right),~
\Gamma=\left(\begin{array}{ccc}0.0041&0.0054&0.0074\\0.0054&0.0067&0.0079\\
    0.0074&0.0079&0.0048\end{array}\right), ~x_0=(1,1,1)^T, ~p_0=(7,7,8)^T.
\end{eqnarray*}
\end{example}

Suppose that the initial liability equity ratio is $l_0/e_0=25$ and  the required leverage ratio is $\rho_1=18$.   One can easily calculate that the initial liability is $l_0=21.1538$, and  the initial equity is $e_0=0.8462$, thus satisfying Assumption~\ref{asmp1}.

We can apply Propositions~\ref{prop-optimality}-\ref{prop-priority} to conclude that the investor should sell the assets to exactly meet the leverage ratio $\rho_1$, and prioritize selling asset 1 than asset 2 during the deleveraging trading. Note that all the elements in the impact matrices are positive in this example. Selling any asset in the portfolio would decrease the prices of all the assets, and thus decrease the investor's final equity. Therefore, the investor would sell assets to precisely meet the liquidity requirement. In addition, compared with the second asset, the first asset is more liquid with lower self price impacts, and less correlated with the third asset with lower averaged cross-impacts. Therefore, the investor would sell more asset 1 than asset 2 to reduce unfavorable price movements, and thus the decrease of his final equity due to his deleveraging trading.

In fact, by solving this example by  Algorithm~\ref{SCOBBA} in Section~\ref{sect:GOA}, one can obtain the global optimal solution
$y^*=(-0.7842,-0.0001,-0.0943)^T$ with the optimal value $e_1(y^*)=0.8287$. It is easy to testify that the value of the leverage constraint at the optimal solution equals $0.000005$, and that $y_{1}^{*}< y_{2}^{*}$. Thus, the leverage constraint of this example is indeed active at its optimal solution, and the first asset is prioritized for selling than asset 2.

\begin{example}\label{exam2}{\rm ~~Consider an instance of Problem~(\ref{OPDP}) with $m=3$ and the following parameters}:
\begin{eqnarray*}
 \Lambda=\left(\begin{array}{rrc}0.0012 &-0.0040& 0.0046\\0.0130& 0.0079& 0.0083\\
-0.0022 &-0.0023 &0.0081 \end{array}\right), ~\Gamma=\left(\begin{array}{rcr} 0.0060 &0.0108 &-0.0062\\
0.0036& 0.0084& 0.0014\\-0.0096& 0.0058& 0.0090\end{array}\right),~x_0=(1,1,1)^T, ~p_0=(7,7,4)^T.
\end{eqnarray*}
\end{example}
Suppose that the initial liability equity ratio is $l_0/e_0=25$ and  the required leverage ratio is $\rho_1=12$.   One can easily calculate that the initial liability is $l_0=17.31$ and  the initial equity is $e_0=0.69$, thus satisfying Assumption~\ref{asmp1}.

Note that $\Gamma$  and $\Lambda$ are asymmetric matrices with negative off-diagonal elements. The trader could increase the price of asset 1 by selling the third asset significantly, due to the fact of  $\lambda_{31} <0$ and $\gamma_{31}<0$. The condition  in Proposition~\ref{prop-optimality} does not hold. However, selling asset 3 could also decrease the price of asset 3 itself and that of asset 2. The trader might still choose not to sell more than the leverage requirement. In fact, by using Algorithm~\ref{SCOBBA} to solve this problem, we obtain its global optimal solution
$y^*=(-1,-0.2241 ,-0.1548)^T$ with the optimal value $e_1(y^*)=0.6855$,  and the value of the leverage constraint  at the optimal solution equals $10^{-6}$. That is to say, the leverage ratio is active at the optimal solution. In addition, it can be also seen from the optimal solution $y^*$ that $y_1^*=-1< -0.2241= y^*_2$. In fact, one can easily check that the four conditions in Proposition~\ref{prop-priority} are satisfied for assets 1 and 2. The first asset is more liquid and less positively correlated with the third asset than asset 2. As a result, the first asset is prioritized for selling than the second asset.

\subsection{Reformulation via Simultaneous Diagonalizability}\label{sec:formulation-SD}
We observe that the two matrices $\Lambda-\frac{1}{2}\Gamma$ and $\Lambda+\frac{1}{2}\Gamma$, respectively, in the objective function and the constraint of Problem (\ref{OPDP})  may be asymmetric, and  indefinite even if they are symmetric. This makes Problem \eqref{OPDP} difficult to be solved.  To remedy this, in this subsection, we present an equivalent reformulation for Problem (\ref{OPDP}) by making use of the spectral decomposition and
the simultaneous diagonalizability of two semi-definite matrices. The idea of simultaneous diagonalizability is also
employed to solve the QCQP problems in \cite{Ben-Tal14,Jiang2018}.

By means of  the symmetric matrices $\hat\Gamma$ and $\hat\Lambda$ defined in \eqref{eqn:sys-matrices}, we first reformulate Problem (\ref{OPDP}) into the following form:
\begin{equation}\label{OPDP-sym}
\begin{split}
\min\limits_{y\in{\mathbb R}^m}~& f(y):=y^T\left(\hat\Lambda-\frac{1}{2}\hat\Gamma\right)y-x_0^T\Gamma y\\
s.~t.~& %g(t):=l_1(y)+\rho_1f(y)\le0,\\
 g(y):=y^T\left(\hat\Lambda+\frac{1}{2}\hat\Gamma+\rho_1\left(\hat\Lambda-\frac{1}{2}\hat\Gamma\right)\right)y+(p_0-\rho_1\Gamma x_0)^Ty+l_0-\rho_1e_0\le0,\\
& -x_0\le y\le0,
\end{split}\tag{$P_1$}
\end{equation}
where $\hat\Gamma+\frac{1}{2}\hat\Lambda$ and $\hat\Gamma-\frac{1}{2}\hat\Lambda$ are symmetric but not necessarily positive definite or semidefinite.
Clearly, Problems~(\ref{OPDP}) and (\ref{OPDP-sym}) are equivalent in the sense that they have the same optimal solution.
When $\hat\Gamma+\frac{1}{2}\hat\Lambda$ and $\hat\Gamma-\frac{1}{2}\hat\Lambda$ are  positive semidefinite, Problem (\ref{OPDP-sym}) is a convex quadratic programming problem which is polynomially solvable in its second-order cone program reformulation. Generally, Problem (\ref{OPDP-sym}) is still a non-convex quadratic program with one quadratic and box constraints, which is NP-hard. 
%\chen{Due to the non-convex quadratic constraint in Problem (\ref{OPDP-sym}), the existing global algorithms in \cite{lblp19,Burer08,Burer09,Burer2012} cannot be directly applied  to solve Problem (\ref{OPDP-sym}).}

We next reformulate Problem (\ref{OPDP-sym}) as an equivalent D.C.  program and discuss how to recover a global solution to Problem (\ref{OPDP-sym}) from the solution of the reformulated problem. Denote by $\nu_i$, $i=1,\ldots, m$,  the  eigenvalues of the matrix $\hat\Lambda-\frac{1}{2}\hat\Gamma$,
% listed in increasing order:
where  $\nu_i<0$ for $i=1,\ldots,s$ and $\nu_i\ge0$ for $i=s+1,\ldots,m$, and by $\eta_i$,   $i=1,\ldots, m$, the corresponding orthogonal unit eigenvectors. We also denote the eigenvalues of the matrix $\hat\Lambda+\frac{1}{2}\hat\Gamma$ by $\mu_i$, $i=1,\ldots, m$, with $\mu_i<0$ for $i=1,\ldots,q$, and $\mu_i\ge0$ for $i=q+1,\ldots,m$, and the corresponding orthogonal unit eigenvectors by $\zeta_i$,   $i=1,\ldots, m$.
%Then we  can separate the matrices to be a minus of two positive semi-definite matrices as follows:
By the spectral decomposition, we can obtain
\begin{eqnarray}\label{B-decomp}
\left\{\begin{array}{l}
\hat\Lambda-\frac{1}{2}\hat\Gamma = B^+-B^-,\quad B^{+}~:=~\sum\limits_{i=s+1}^m\nu_i\eta_i\eta_i^T ,\quad B^{-}~:=~\sum\limits_{i=1}^s-\nu_i\eta_i\eta_i^T , \\
\hat\Lambda+\frac{1}{2}\hat\Gamma = A^+-A^-,\quad A^{+}~:=~\sum\limits_{i=q+1}^m\mu_i\zeta_i\zeta_i^T ,\quad A^{-}~:=~\sum\limits_{i=1}^q-\mu_i\zeta_i\zeta_i^T .
\end{array}\right.%\label{A-decomp}
\end{eqnarray}
%According to the well-known result regarding the simultaneous diagonalizability of two semi-definite matrices from \cite{Newcomb61} (see Lemma~\ref{lem-SD} in the e-companion to this paper),
Recall a well-known result  regarding the simultaneous diagonalizability of two semi-definite matrices from \cite{Newcomb61}.
\begin{lemma}\label{lem-SD}
Let $A$ and $B$ be two $n\times n$ real symmetric positive semi-definite matrices. Then $A$ and $B$ are simultaneously diagonalizable, i.e., there exists a nonsingular matrix $D$ such that both $D^T AD$ and $D^T BD$ are diagonal.
\end{lemma}
According to Lemma~\ref{lem-SD}, because $B^-$ and $A^-$ are symmetric positive semi-definite matrices,  $B^-$ and $A^-$ are  simultaneously
diagonalizable. In particular, we can identify a non-singular matrix $D$ to simultaneously diagonalize these two matrices as
\begin{eqnarray}\label{SD}
 D^T B^-D ={\rm diag}(\delta_1,\ldots,\delta_s,0,\ldots,0),\quad\quad D^T A^-D ={\rm diag}(\theta_1,\ldots,\theta_r,0,\ldots,0),
\end{eqnarray}
where $r={\rm rank} (B^-+A^-)$ with $\max\{s,q\}<r\le s+q$, $0<\delta_i\le 1$ for $i=1,\ldots,s$,  $\theta_i=1-\delta_i$ for $i=1,\ldots,s$  and $\theta_{i}=1$ for  $i=s+1,\ldots,r$, and there are $q$ nonzero numbers in $\theta_1,\ldots,\theta_r$ (see the proof of Lemma~\ref{lem-SD} in the e-companion to this paper).
Using the transformation $y=Dz$, we then obtain the following proposition.
\begin{proposition}\label{prop-SD}
Problem (\ref{OPDP-sym}) %can be equivalently  transformed into the following form
has the same optimal value as the following  D.C.  program
\begin{equation}\label{OPDP-SD}
\begin{split}
\min\limits_{z\in{\mathbb R}^m}~& \hat{f}(z):= z^T(D^TB^+D)z-x_0^T\Gamma Dz-\sum_{i=1}^{s}\delta_iz_i^2\\
{\rm s.~t.}~& \hat g(z):= \psi(z)-\sum_{i=1}^{r}\theta_iz_i^2-\rho_1\sum_{i=1}^{s}\delta_iz_i^2\le0,\\
 & z\in{\cal Z}:=\left\{ z\in{\mathbb R}^m\mid  -x_0\le Dz\le0 \right\},
\end{split}\tag{$\hat{P}$}
\end{equation}
where $\psi(z)$ is a quadratic convex function given by
\begin{eqnarray}\label{defn-psi-fun}
\psi(z):= z^TD^T\left(A^++\rho_1B^+\right)Dz+(p_0-\rho_1\Gamma x_0)^TDz+l_0-\rho_1e_0.
\end{eqnarray}
Furthermore, if $z^*$ is a global optimal solution to Problem (\ref{OPDP-SD}), then $y^*=Dz^*$ is a global optimal solution to Problem~(\ref{OPDP-sym}).
\end{proposition}

%\chen{
%\begin{remark}\label{lem-SD}
% If $\hat\Gamma+\frac{1}{2}\hat\Lambda$ is positive semidefinite, then Problem (\ref{OPDP-sym}) can be equivalently  transformed into the following form
%\begin{equation}\label{OPDP-DC}
%\begin{split}
%\min\limits_{y\in{\mathbb R}^m,t\in{\mathbb R}^s}~&  y^TB^+y-x_0^T\Gamma y-\|z\|_2^2\\
%{\rm s.~t.}~& \varphi(y) -\rho_1\|z\|_2^2\le0,\\
% & z=Cy,~~-x_0\le y\le0 ,
%\end{split}\tag{$\tilde{P}$}
%\end{equation}
%where $C=\left(\sqrt{-\nu_1}\eta_1,\ldots,\sqrt{-\nu_s}\eta_s\right)^T$, and $\varphi(y)$ is a quadratic convex function given by
%\begin{eqnarray*}\label{defn-varphi-fun}
%\varphi(y):= y^T\left(\hat\Gamma+\frac{1}{2}\hat\Lambda+\rho_1B^+\right)y+(p_0-\rho_1\Gamma x_0)^Ty+l_0-\rho_1e_0.
%\end{eqnarray*}
%\end{remark}
%}

It is worth mentioning  that Problem (\ref{OPDP-SD}) is  still non-convex  and its global optimal solution is hard to obtain. However,  the objective function and the constraint  of Problem (\ref{OPDP-SD}) have simple separable concave quadratic terms $-z_i^2$. Based on this special structure, we will, in the subsequent sections, develop an effective algorithm  to find a global optimal solution to Problem (\ref{OPDP-SD}).

\section{The SCO Method and Its Convergence}\label{sec:SCO}
In this section, we propose a successive convex optimization (SCO) algorithm for Problem \eqref{OPDP-SD}, and prove the convergence of the SCO algorithm to a KKT point of Problem (\ref{OPDP-SD}).

\subsection{Quadratic Convex Approximation}\label{subsec:QCA}
 We first present a  quadratic convex approximation for Problem (\ref{OPDP-SD}) by linearizing the concave quadratic terms in
the objective function and the constraint. This idea of linearizing the concave terms in the D.C. program has been used in the literature (cf. \cite{Hong11,LeThi14,lblp19,ldpjl20}).

 Let  $z_l,z_u\in{\mathbb R}^{m}$  denote, respectively, the lower and upper bounds of $z=D^{-1}y$ for $y\in[-x_0,0]$  given by
\begin{eqnarray}\label{z-bounds}
z_l^i=-\sum\limits_{j:\hat{D}_{ij}>0}\hat{D}_{ij}x_{0,j}, \quad  z_u^i=-\sum\limits_{j:\hat{D}_{ij}<0}\hat{D}_{ij}x_{0,j}, \quad i=1,\cdots,m,
\end{eqnarray}
where $\hat{D}_{ij}$ denotes the element of the matrix $D^{-1}$.
For an arbitrary given $\xi_i\in [z_l^i,z_u^i]$, one has
%\chen{we construct a linear term of $z_i$, which serves as an upper bound of the quadratic term $ -z_i^2$, with the following induction:}
\begin{equation}\label{linearization}
 -z_i^2\le -\xi_i^2-2\xi_i(z_i-\xi_i)=-2\xi_iz_i+\xi^2_i,\quad\forall z_i\in[z_l^i,z_u^i].
 \end{equation}
Using the  linear term $(-2\xi_iz_i+\xi^2_i)$  to approximate the concave quadratic terms $-z^2_i$ in the objective function and the constraint  of Problem (\ref{OPDP-SD}), we have the following quadratic convex approximation,
    \begin{eqnarray}\label{ODP-linearappro}
  &\min\limits_{z\in{\mathbb R}^m}& \hat{f}_\xi(z):= z^T(D^TB^+D)z-x_0^T\Gamma Dz +\sum_{i=1}^{s}\delta_i\left(-2\xi_iz_i+\xi^2_i\right)\\\nonumber
         &s.~t.&  \hat g_{\xi}(z):=\psi(z)+\sum_{i=1}^{r}\theta_i\left(-2\xi_iz_i+\xi^2_i\right)+\rho_1\sum_{i=1}^{s}\delta_i\left(-2\xi_iz_i+\xi^2_i\right)\le 0, \\\nonumber
         && z\in{\cal Z},
\end{eqnarray}
 where $\xi=(\xi_1,\ldots,\xi_r)^T$ with $\xi_i\in [z_l^i,z_u^i]$, $i=1,\ldots,r$,  and $\psi(z)$ is a convex function defined in (\ref{defn-psi-fun}).
 Note that inequality (\ref{linearization})  implies that $\hat g_{\xi}(z)\ge \hat g(z)$,
thus the optimal solution to problem (\ref{ODP-linearappro}) is also feasible to Problem (\ref{OPDP-SD}). Hence, the objective function value $\hat f(z)$ at the optimal solution to Problem (\ref{ODP-linearappro})  provides an upper bound for Problem (\ref{OPDP-SD}).

We denote the feasible sets of  Problems~(\ref{OPDP-SD}) and \eqref{ODP-linearappro} and the interior of \eqref{ODP-linearappro} by
\begin{eqnarray*}
 \hat{\cal F}=\{z\in{\cal Z}\mid \hat{g}(z)\le0\},\quad \hat{\mathcal{F}}_{\xi}=\left\{ z\in{\cal Z}\mid   \hat{g}_{\xi}(z)\le 0\right\},\quad\text{int}\hat{\mathcal{F}}_{\xi}=\left\{ z\in{\cal Z}\mid  \hat{g}_{\xi}(z)< 0 \right\}.
\end{eqnarray*}
We next present a simple property of the sets $\hat{\mathcal{F}}_{\xi}$ and ${\rm int}\hat{\mathcal{F}}_{\xi}$.
\begin{lemma}\label{lem-F-nonemptiness}
 Let   $\bar z\in\hat{\cal F}$ and $\bar\xi_i=\bar z_i$ for $i=1,\ldots,r$. We have the following.
\begin{description}
 \item[(i)]  $\hat{\cal F}_{\bar\xi}$ is a nonempty closed convex set in $\hat{\cal F}$.
 \item[(ii)] ${\rm int}\hat{\cal F}_{\bar\xi}\neq\emptyset$.
 \end{description}
\end{lemma}

From Lemma~\ref{lem-F-nonemptiness}, we see that  Problem~(\ref{ODP-linearappro}) is feasible and well-defined
and that the Slater constraint qualification holds for Problem~(\ref{ODP-linearappro}).  The following proposition follows immediately.
  \begin{proposition}\label{prop2-optimality}
 Let $\bar z\in\hat{\cal F}$ and $\bar\xi_i=\bar z_i$ for $i=1,\ldots,r$.
 If $\bar z$ is an optimal solution to  Problem (\ref{ODP-linearappro}) with  $\xi=\bar\xi$,
 then  $\bar{z}$  is a  KKT point   of  Problem (\ref{OPDP-SD}).
 \end{proposition}

\subsection{The SCO Algorithm}\label{subsec:SCO}
%\chen{We thus propose our SCO algorithm based on the quadratic convex approximation presented in the above subsection. In particular, with an initial feasible solution $z^{0}$, we update the vectors $z^{k}$ and $\xi^k$ to approximate the solution to \eqref{OPDP-SD}. For the $k$-th iteration, we iterate over the following two steps: first, update the new vector $\xi^{k+1}$ as $\xi^{k+1}_{i} = z^{k+1}_{i}$, $i=1,\cdots,r$; then, solve the quadratic convex problem \eqref{ODP-linearappro} with $\xi = \xi^k$ and document the solution as $z^{k+1}$.} We  describe  the  SCO algorithm  for  problem~(\ref{OPDP-SD}) in Algorithm \ref{SCO}.

 Next, we  describe  the  SCO algorithm  for  Problem~(\ref{OPDP-SD}) based on the convex approximation (\ref{ODP-linearappro}).
\begin{algorithm} [!ht]\caption{\bf [SCO Algorithm $(SCO(z^0,\epsilon))$]}
\label{SCO} {~}\vskip 0pt {\rm
\begin{description}
\item[\bf Input: ]   Initial feasible point $z^0\in\hat{\cal F}$ and  stopping criterion  $\epsilon>0$. % Initialize $k=0$.
\item[Step  0]  Set  $\xi_i^0=z_i^0$ for $i=1,\ldots,r$ and $\xi^0=(\xi_1^0,\ldots,\xi_r^0)^T$.  Set $k=0$.
\item[Step  1] Solve problem~(\ref{ODP-linearappro}) with  $\xi=\xi^k$ to get the optimal solution $z^{k+1}$. Set $\xi^{k+1}=(\xi_1^{k+1},\ldots,\xi_r^{k+1})^T$ with $\xi_i^{k+1}=z_i^{k+1}$ for $i=1,\ldots,r$.
  \item[Step 2] If $\|\xi^{k+1}-\xi^k\|> \epsilon$,  then set $k=k+1$ and go back to step 1;
   Otherwise, stop and output $z^{k+1}$ as the final output.
\end{description}
 }
\end{algorithm}

It should be pointed out that the SCO algorithm needs an initial solution $z^0\in\hat{\cal F}$. By Proposition~\ref{prop:property-problem}, $-x_0$ is a feasible solution of Problem (\ref{OPDP}). Note that Problems (\ref{OPDP}) and (\ref{OPDP-sym}) have the same feasible set, and that $y=Dz$ with $D$ being non-singular. Thus, $z^0=D^{-1}(-x_0)$ is also a feasible solution to Problem (\ref{OPDP-SD}).

%Let $\{z^k\}$ be the sequence generated by Algorithm \ref{SCO}. Note that since $\{Pz^k\}\subseteq[-x_0,0]$ and $P$ is nonsingular, the sequence $\{z^k\}$ is bounded and hence it  has at least one accumulation point.
We show in the following lemma that the sequence $\{z^k \}$ generated by Algorithm \ref{SCO} always lie inside the feasible region of Problem (\ref{OPDP-SD}), with the corresponding objective value improved by every iteration.
\begin{lemma}\label{lem1-SCA}
The sequences $\{z^k \}$ and $\{\xi^k \}$ generated by Algorithm \ref{SCO} satisfy that
\begin{description}
\item[(i)] $\{z^k \}\subseteq\hat{\cal F}$. %$z^k$ is always feasible to Problem \eqref{ODP-linearappro}.
\item[(ii)] $\hat{f}(z^{k})-\hat{f}(z^{k+1}) \ge \sum\limits_{i=1}^{s}\delta_i(\xi^{k+1}_i-\xi^k_i)^2.$
\end{description}
\end{lemma}

Note that, because $\{Dz^k\}\subseteq[-x_0,0]$ with the non-singular matrix $D$, the sequence $\{z^k\}$ is bounded, and thus has at least one accumulation point. Based on Lemma~\ref{lem1-SCA}, we obtain the following lemma.
\begin{lemma}\label{lem3-SCA}
 Let the   sequence   $\{(z^k,\xi^k) \}$ be generated by Algorithm \ref{SCO} with an accumulation point
  $(\hat{z},\hat{\xi})$.  Then, (i)  $\hat{z}\in {\cal F}_{\hat\xi}$ and $\hat{\xi}_i=\hat{z}_i$ for $i=1,\ldots,r$, and (ii) $\hat{f}_{\hat\xi}(z)\ge \hat{f}_{\hat\xi}(\hat{z})$ for any $z\in {\cal F}_{\hat\xi}$.
\end{lemma}

Lemma~\ref{lem3-SCA} indicates that the accumulation point $(\hat{z},\hat{\xi})$ of the sequence $\{(z^k,\xi^k)\}$ generated by Algorithm~\ref{SCO}
satisfies that $\hat{z}$ is the optimal solution to Problem (\ref{ODP-linearappro}) with $\xi=\hat{\xi}$.
  Combining  Lemma~\ref{lem3-SCA} with Proposition~\ref{prop2-optimality}, we immediately obtain  the following convergence result for Algorithm~\ref{SCO}.
\begin{theorem}\label{thm-SCO}
 Let $\bar{z}$ be an accumulation point of the sequence $\{z^k\}$  generated by Algorithm \ref{SCO} with $\epsilon=0$. % and  $\bar{y}=P\bar z$.
 Then $\bar{z}$ is a  KKT point   of  Problem (\ref{OPDP-SD}).
 \end{theorem}

Theorem \ref{thm-SCO} shows  that  the SCO algorithm converges to a KKT point of  Problem (\ref{OPDP-SD}).
In the next section, we would utilize this property of the SCO algorithm to design a global algorithm for Problem (\ref{OPDP-SD}). In particular, the solution derived by the SCO algorithm can be used as an upper bound in our branch-and-bound approach to identify the global optimizer. In addition, we will see in our numerical experiments in Section~\ref{sec:Experiment} that the KKT point identified by the SCO algorithm for Problem \eqref{OPDP-SD}  always provides a good approximation for the global optimal solution.

\section{A Global Optimization Algorithm }\label{sect:GOA}
In the section,  we develop an algorithm that identifies the global optimal solution to Problem (\ref{OPDP-SD}) within a pre-specified $\epsilon$-tolerance by integrating the SCO algorithm with a simple convex relaxation  and branch-and-bound framework.
We  also establish the convergence of the algorithm and estimate its complexity.

\subsection{The Quadratic Convex Relaxation}\label{sec:QCR}

In this subsection, we present a simple quadratic convex relaxation for Problem~(\ref{OPDP-SD}) via the convex envelope technique, and then estimate the gap between   it and the original problem.

To start, we  consider a restricted version of Problem~(\ref{OPDP-SD}),
  where the variables $z_i$ $(i=1,\ldots,r)$ are in a sub-rectangle $[l,u]$ with $l,u\in {\mathbb R}^{r}$:
 \begin{eqnarray}\label{OPD:Bounded}
 &\min\limits_{z\in{\mathbb R}^m}& \hat{f}(z):= z^T(D^TB^+D)z-x_0^T\Gamma Dz-\sum_{i=1}^{s}\delta_iz_i^2\\\nonumber
&{\rm s.~t.}& \hat g(z):= \psi(z)-\sum_{i=1}^{r}\theta_iz_i^2-\rho_1\sum_{i=1}^{s}\delta_iz_i^2\le0,\nonumber\\
 && -x_0\le Dz\le0,\quad z_i\in [l_i,u_i],~i=1,\ldots,r,\nonumber
\end{eqnarray}
where $[l_i,u_i]\subseteq[z_l^i,z_u^i]$  and $z_l^i,z_u^i$  are given in (\ref{z-bounds}).  Let $t_i=z_i^2$ for $ i=1,\ldots,r$. According to \cite{McCormick76}, the  convex envelope  of $t_i=z_i^2$ on $[l_i,u_i]$ is $\left\{(t_i,z_i):z_i^2\le t_i,~t_i \le (l_i+u_i)z_i-l_iu_i\right\}$.
Replacing the negative quadratic term $-z_i^2$ in the objective function and the constraint of Problem (\ref{OPD:Bounded}) with its convex envelope,
we can derive the following quadratic convex relaxation for Problem~(\ref{OPD:Bounded}):
\begin{eqnarray}\label{ODP:Bounded-CR}
&\min\limits_{z\in{\mathbb R}^m,t\in{\mathbb R}^r}  & z^T(D^TB^+D)z-x_0^T\Gamma Dz-\sum_{i=1}^{s}\delta_it_i , \\\nonumber
&{\rm s.~t.}&  \psi(z)-\sum_{i=1}^{r}\theta_it_i-\rho_1\sum_{i=1}^{s}\delta_it_i\le0,\\
 && -x_0\le Dz\le0,\quad z_i\in [l_i,u_i],~i=1,\ldots,r,\nonumber\\
&& z_i^2\le t_i, \quad   t_i \le (l_i+u_i)z_i-l_iu_i,\quad i=1,\ldots,r. \nonumber
 \end{eqnarray}
It is worth mentioning that the last constraints on $(t_i,z_i)$ in (\ref{ODP:Bounded-CR}) are so-called secant cuts originally introduced in \cite{Saxena11}, and that the above convex envelope technique  was also used in \cite{lblp19,ldpjl20}.

 As pointed out in \cite{lblp19}, although there exist many other strong relaxation models for QCQP,  they usually involve more intensive computation.
In this work, we  combine the relaxation model (\ref{ODP:Bounded-CR}) with other simple optimization techniques to develop a global algorithm for Problem~(\ref{OPDP-SD}).
Our choice is based on the special structure of the relaxation model (\ref{ODP:Bounded-CR}) that has the separability of the constraints
on $(t_i,z_i)$. Such a structure  allows us to adopt the adaptive branch-and-cut technique based on partition on the variables $z_1,\ldots,z_r$ in the design of the global algorithm. Moreover, the relaxation (\ref{ODP:Bounded-CR}) has the good approximation behavior as shown in our next theorem, which
compares the optimal values of Problem~\reff{OPD:Bounded} and its relaxation problem~(\ref{ODP:Bounded-CR}).

 \begin{theorem}\label{thm1-CR}
Let $\hat{f}^{*}_{[l,u]}$ and $\hat{v}^*_{[l,u]}$ be  the optimal value  of Problem~(\ref{OPD:Bounded}) and its relaxation~(\ref{ODP:Bounded-CR}), respectively.  Let $(\bar{z},\bar{t})$ be the optimal solution to Problem (\ref{ODP:Bounded-CR}).  Then,
\begin{eqnarray*}
&& \hat{f}(\bar z)-\hat{f}^*_{[l,u]}\le \hat{f}(\bar z)-\hat{v}^*_{[l,u]}\le\sum_{i=1 }^{s}(\bar t_i-  \bar z_i^2)
\le \frac{s}{4}\|u-l\|^2_\infty,\\
&& \hat{g}(\bar z)\le (r+\rho_1s)\max_{i=1,\ldots,r}\{\bar t_i-  \bar z_i^2\}\le\frac{r+\rho_1s}{4}\|u-l\|^2_\infty,
 \end{eqnarray*}
 where $\|\cdot\|_\infty$ denotes the $\ell_\infty$-norm on ${\mathbb R}^{r}$  defined by
  $\|\xi\|_\infty=\max\limits_{i=1,\ldots,r}|\xi_i|$.
\end{theorem}
%\proof %Since problem  (\ref{ODP:Bounded-CR}) is a convex relaxation of problem (\ref{OPD:Bounded}),
%We first have $\hat{f}^{*}_{[l,u]}\ge \hat{v}^*_{[l,u]}(\bar{y})$. Note that $0<\delta_i\le1$ for $i=1,\ldots,s$.
%  It follows  that
%\begin{eqnarray*} \nonumber
%  \hat{f}(\bar z)-\hat{f}^*_{[l,u]}&\le& \hat{f}(\bar z)-\hat{v}^*_{[l,u]}=\sum_{i=1 }^{s}\delta_i(\bar t_i-  \bar z_i^2) \\
%    & \le&\sum_{i=1 }^{s} (\bar t_i-  \bar z_i^2) \label{thm2-CR-ineq}\\\nonumber
%  &\le& \sum_{i=1 }^{s} [-\bar z_i^2+(l_i+u_i)\bar z_i -l_iu_i] \\\nonumber
% &\le& \frac{1}{4}\sum_{i=1 }^{s}(u_i-l_i)^2\le \frac{s}{4}\|u-l\|^2_\infty,
%  \end{eqnarray*}
% where the second and third inequalities follow  from the constraints on $(t_i, z_i)$ in (\ref{ODP:Bounded-CR}).
%
% Similarly, using the feasibility of $(\bar{z},\bar{t})$ and the fact that $0<\theta_i\le1$ and $0<\delta_i\le1$ for each $i$, we can follow that
%\begin{eqnarray*} \nonumber
%\hat{g}(\bar z)&=&\psi(\bar z)-\sum_{i=1}^{r}\theta_i\bar t_i-\rho_1\sum_{i=1}^{s}\delta_i\bar t_i+\sum_{i=1}^{r}\theta_i(\bar t_i-\bar z_i^2) +\rho_1\sum_{i=1}^{s}\delta_i(\bar t_i-\bar z_i^2)\\
%&\le& \sum_{i=1}^{r}(\bar t_i-\bar z_i^2) +\rho_1\sum_{i=1}^{s}(\bar t_i-\bar z_i^2)\\
%&\le&(r+\rho_1s)\max_{i=1,\ldots,r}\{\bar t_i-  \bar z_i^2\}\le\frac{r+\rho_1s}{4}\|u-l\|^2_\infty.
%  \end{eqnarray*}
%This completes the proof.  \qed

%
To define how good an approximate solution is, we extend the $\epsilon$-solution introduced in \cite{Lor82} to the following definition.
\begin{definition}\label{defn-local-minimizer}
Let $\epsilon>0$ and $\hat{f}^{*}$  be  the optimal value  of Problem~(\ref{OPDP-SD}).
The point $\bar z\in {\cal Z}$ is said to be an $\epsilon$-optimal solution  to Problem (\ref{OPDP-SD}) if $\hat{g}(\bar z)\le\epsilon$ and
$\hat{f}(\bar z)-\hat{f}^{*}\le\epsilon$.
\end{definition}

%We remark that the above definition can be viewed as an extension of the so-called $\epsilon$-solution introduced in \cite{Lor82}.
From  Theorem~\ref{thm1-CR}, we have the following corollary immediately.
\begin{corollary}\label{cor-CR}
Let $ (\bar{z},\bar{t}) $ be the optimal solution to Problem (\ref{ODP:Bounded-CR}) and
     $\epsilon>0$. If  $\max\limits_{i=1,\ldots,r}\{\bar t_i-  \bar z_i^2\}\le\frac{\epsilon}{r+\rho_1s}$,
     then $\bar{z}$ is an $\epsilon$-optimal solution of Problem~(\ref{OPD:Bounded}).
 \end{corollary}

Theorem~\ref{thm1-CR} indicates that %when $\|u-l\|_\infty$ is very small,
when the length of the longest edge of rectangle $[l, u]$ is sufficiently short, %less than or equal to $2\sqrt{\epsilon/r}$,
  the relaxed model (\ref{ODP:Bounded-CR}) can provide a good approximation  to %the original
problem (\ref{OPD:Bounded}). Moreover, from  Theorem  \ref{thm1-CR}, we see that if
$\|u-l\|_\infty \le \frac{2\sqrt{\epsilon}}{\sqrt{r+\rho_1s}}$,
then $\hat{g}(\bar z)\le\epsilon$ and $\hat{f}(\bar{z})-\hat{f}^{*}_{[l,u]}\le\epsilon$. Thus, $\bar{z}$ can be viewed as an $\epsilon$-optimal solution to problem (\ref{OPD:Bounded}).
%Motivated by this observation, we can divide the initial interval $[z_l^i,z_u^i]$ into
%$\lceil{\frac{(z_u^i-z_l^i)\sqrt{r+\rho_1s}}{2\sqrt{\epsilon}}}\rceil$ subintervals such that each subinterval has a width of $\frac{2\sqrt{\epsilon}}{\sqrt{r+\rho_1s}}$, where $z_l^i,z_u^i$  are given in (\ref{z-bounds}).
%Then we can solve the  relaxed problem~(\ref{ODP:Bounded-CR}) for every sub-rectangle  to obtain an $\epsilon$-optimal solution for  the restricted
%problem (\ref{OPD:Bounded}) over this sub-rectangle. After that, we can choose the best solution from the obtained approximate solutions as the final solution,
%which is clearly an $\epsilon$-optimal solution to Problem (\ref{OPDP-SD}).
%From the above discussion, we can see that such partitioning procedure can find an $\epsilon$-optimal solution to Problem (\ref{OPDP-SD}) in
%$\mathcal{O} \left( N\prod_{i=1}^{r}\lceil\frac{\sqrt{r+\rho_1s}(z_u^i-z_l^i)}{2\sqrt{\epsilon}}\rceil\right)$
% time, where $N$ is the complexity to solve problem (\ref{ODP:Bounded-CR}). However, such a partitioning procedure is not effective.
%In next subsection, based on Theorem~\ref{thm1-CR}, we shall discuss how to integrate the SCO approach  with  branch-and-bound techniques to develop an effective global algorithm for Problem (\ref{OPDP-SD}).

\subsection{The SCOBB Algorithm}\label{sect_IA}

In this subsection, we present a global algorithm (called SCOBB) for Problem (\ref{OPDP-SD}) that integrates  the SCO approach with the B\&B framework based on  the quadratic convex relaxation (\ref{ODP:Bounded-CR}) and the adaptive branch-and-cut rule from \cite{lblp19}.  The SCOBB algorithm for   Problem~(\ref{OPDP-SD}) is described in Algorithm~\ref{SCOBBA}. \\[1pt]

\begin{breakablealgorithm}\caption{\bf (The SCOBB Algorithm)}
\label{SCOBBA} {~}\vskip 0pt {\rm
\begin{description}
\item[\bf Input:] $\Lambda,\Gamma, p_0,x_0,\rho_1$,   and stopping criteria $\epsilon>0$.
\item[\bf Output:] $\epsilon$-optimal solution $z^*$.
\item [\bf Step 0]  ({\bf Initialization})
\begin{description}
 \item [(i)] Compute positive semi-definte matrices $B^-$ and $A^-$ by (\ref{B-decomp}).

 \item [(ii)] Compute a nonsingular matrix $D$ such that $D^TB^-D$ and $D^TA^-D$ are diagonal.
Solve the equation $Dz=-x_0$ to obtain a solution $\bar z$.

 \item [(iii)] Let $l^0=(z_l^1,\ldots, z_l^r)^T$, $u^0=(z_u^1,\ldots, z_u^r)^T$,  where $z_l^i$ and $z_u^i$  are computed by (\ref{z-bounds}).
\end{description}
\item[\bf Step 1] Find a KKT point  $z^*$ of Problem (\ref{OPDP-SD}) by running SCO$(z^0,\epsilon)$ with $z^0=\bar{z}$. Set the first upper bound $v^*=\hat{f}(z^*)$.

 \item [\bf Step 2]  Solve  problem~(\ref{ODP:Bounded-CR}) over $[l=l^0,u=u^0]$ to obtain the optimal solution $(z^0,t^0)$ and the first lower bound $v^0$.
      %Set $y^0=Pz^0$.
      If $\hat{g}(z^0)\le \epsilon$   and $\hat{f}(z^0)< v^*$, then update  upper bound $v^*=\hat{f}(z^0)$ and solution $z^*=z^0$.
        Set $k=0$, $\Delta^k:=[l^k,u^k]$, $\Omega:=\left\{[\Delta^k, v^k, (z^k,t^k)]\right\}$.
\item [\bf Step 3] {\bf While} $\Omega\neq\emptyset$  {\bf Do} (the main loop)
\begin{description}
\item[\bf (S3.1)] ({\bf Node Selection}) Choose a node $[\Delta^k, v^k, (z^k,t^k)]$ from $\Omega$
with the smallest lower bound $v^k$  and delete it from $\Omega$.

\item[\bf (S3.2)] ({\bf Termination}) If $v^k\ge v^*-\epsilon$, then $z^*$ is an  $\epsilon$-optimal  solution to Problem (\ref{OPDP-SD}), stop.

\item [\bf (S3.3)] ({\bf Partition}) Choose $i^*$ maximizing $t_i^k-(z_i^k)^2$ for $i=1,\ldots,r$.
Set $w_{i^*}=\frac{l^k_{i^*}+u^k_{i^*}}{2}$,
$$\Phi_k(w_{i^*})=\left\{(t_{i^*},z_{i^*})\left |
\begin{array}{l}
t_{i^*}> (l_{i^*}^k+w_{i^*})z_{i^*}-l_{i^*}^kw_{i^*}\\t_{i^*}> (w_{i^*}+u^k_{i^*})z_{i^*}-w_{i^*}u^k_{i^*}
\end{array}\right.\right\}.$$
If $(t^k_{i^*},z^k_{i^*})\in \Phi_k(w_{i^*})$, then set the branching point $\beta_{i^*}=w_{i^*}$; else set $\beta_{i^*}=z^k_{i^*}$. Partition $\Delta^k$ into two sub-rectangles $\Delta^{k_1}$ and
$\Delta^{k_2}$ along the edge $[l^k_{i^*},u^k_{i^*}]$ at point $\beta_{i^*}$.

 \item[\bf (S3.4)] For $j=1,2$, if  problem~(\ref{ODP:Bounded-CR}) over $\Delta^{k_j}$ is feasible, then
 solve problem~(\ref{ODP:Bounded-CR})  over $\Delta^{k_j}$  to get an  optimal solution $(z^{k_j},t^{k_j})$ and optimal value $v^{k_j}$.
Set $\Omega=\Omega\cup\left\{[\Delta^{k_1}, v^{k_1},(z^{k_1},t^{k_1})]\right\}\cup\left\{[\Delta^{k_2}, v^{k_2},(z^{k_2},t^{k_2})]\right\}$.

\item[\bf (S3.5)] ({\bf Restart SCO}) Set $\tilde{z}=\arg\min\{\hat{f}(z^{k_1}), \hat{f}(z^{k_2})\}$.
If $\hat{g}(\tilde{z})\le \epsilon$ and $\hat{f}(\tilde{z})\le v^*-\epsilon$, then  find a  KKT point  $\bar z^k$   of Problem (\ref{OPDP-SD}) by running $SCO(z^0,\epsilon)$ with  $z^0=\tilde{z}$, update solution  $z^*=\arg\min\{\hat{f}(\tilde{z}), \hat{f}(\bar{z}^k)\}$ and upper bound $v^*=\hat{f}(z^*)$.
\item[\bf (S3.6)]  ({\bf Node deletion}) Delete from $\Omega$ all the nodes $[\Delta^j, v^{j},(z^j,t^j)]$ with $v^j\geq v^*-\epsilon$.   Set $k=k+1$.
\item[\bf End while]
\end{description}
\end{description}
}\end{breakablealgorithm}

\vskip 18pt
We remark that there are two main differences between the SCOBB algorithm and the existing global algorithms for QP with non-convex constraints (cf. \cite{Audet00,Linderoth05}):

(i) The SCOBB algorithm utilizes SCO to compute upper bounds for Problem (\ref{OPDP-SD}) to accelerate the convergence, and  improves the upper bound by restarting SCO under certain circumstance;

 (ii) Based on the special structure of the relaxation (\ref{ODP:Bounded-CR}), which is the separability of the constraints
on $(t_i,z_i)$, we adopt the adaptive branch-and-cut rule from \cite{lblp19} to cut off the optimal solution to the relaxation problem corresponding to the node with the smallest lower bound after each iteration, so that the lower bound will be improved (see Figure 1 for an illustration in  \cite{lblp19}).

We next present a technical result  for the sequence $\{(z^k,t^k)\}$ generated by Algorithm~\ref{SCOBBA}.
%Recall that $(z^k,t^k)$ is the optimal solution of problem~(\ref{ODP:Bounded-CR}) over $\Delta^k$. From the proof of Theorem~\ref{thm1-CR}, we  immediately have the following lemma.
\begin{lemma} \label{lemma-SCOBBA}
\begin{description}
\item[\rm (i)] $t_i^k-(z_i^k)^2\le \frac{1}{4}(u_i^k-l_i^k)^2$, $i=1,\ldots,r$ for all $k$.

\item[\rm (ii)]  At the $k$-th iteration, if $\max\limits_{i=1,\ldots,r}\{t_i^k-(z_i^k)^2\}\le \frac{\epsilon}{r+\rho_1s}$, then
Algorithm~\ref{SCOBBA}  stops and $z^*$ is an  $\epsilon$-optimal solution  to Problem (\ref{OPDP-SD}).
\end{description}
\end{lemma}

%%The following lemma can be proved using the similar arguments as for Lemma 3.2 in \cite{lblp19}.
%\begin{lemma}\label{lemma-SCOBBA-stop}
%At the $k$-th iteration, if $\max\limits_{i=1,\ldots,r}\{t_i^k-(z_i^k)^2\}\le \frac{\epsilon}{r+\rho_1s}$, then
%Algorithm~\ref{SCOBBA}  stops and $z^*$ is a global  $\epsilon$-approximate solution  to problem (\ref{OPDP-SD}).
%\end{lemma}
%\proof The proof is similar to the proof of Lemma 3.2 in \cite{lblp19} based on Theorem~\ref{thm1-CR}. \qed

%
%\proof  Since $(z^k,t^k)$ is the optimal solution of problem~(\ref{ODP:Bounded-CR}) over $\Delta^k$ and $v^k$ is its optimal value, by Theorem~\ref{thm1-CR}, we can infer that $\hat{f}(z^k)-v^k\le\epsilon$ and $\hat{g}(z^k)\le\epsilon$. Thus, $z^k$ is $\epsilon$-feasible for problem (\ref{OPDP-SD}).
%From Steps 2 and (S3.5)  of the algorithm, we see that $\hat{f}(z^k)\ge v^*=\hat{f}(z^*)$ for all $k$. Thus,
%\begin{eqnarray}\label{lemma-SCOBBA-ineq1}
%v^k- v^*\ge v^k-\hat{f}(z^k)\ge-\epsilon,
%\end{eqnarray}
% so the stopping criterion of the algorithm is satisfied and then the algorithm   stops. Let $\hat{f}^*$ denotes the optimal value of problem (\ref{OPDP-SD}). By Step (S3.1), $v^k$ is the smallest lower bound. Thus $\hat{f}^*\ge v^k$. It then  follows from (\ref{lemma-SCOBBA-ineq1}) that
% \begin{eqnarray*}\label{lemma-SCOBBA-ineq2}
%    \hat{f}(z^*)=v^*\le v^k+\epsilon\le \hat{f}^*+\epsilon.
%\end{eqnarray*}
%This, together with   $\hat{g}(z^*)\le\epsilon$, implies that   $z^*$ is a global  $\epsilon$-solution to problem (\ref{OPDP-SD}). \qed

%\vskip 6pt
We can now establish  the convergence of Algorithm~\ref{SCOBBA} %to an $\epsilon$-optimal solution of Problem (\ref{OPDP-SD})
based on Lemma~\ref{lemma-SCOBBA}.
 \begin{theorem}\label{thm1-algo1-complexity}
Algorithm~\ref{SCOBBA} can identify an  $\epsilon$-optimal solution to Problem (\ref{OPDP-SD}) via solving at most
$\prod\limits_{i=1}^{r}\left\lceil{\frac{\sqrt{r+\rho_1s}(z_u^i-z_l^i)}{2\sqrt{\epsilon}}}\right\rceil$ relaxed subproblem (\ref{ODP:Bounded-CR}).
\end{theorem}

Theorem~\ref{thm1-algo1-complexity} indicates that Algorithm~\ref{SCOBBA} enjoys a worst-case complexity bound  $\mathcal{O} \left(\prod\limits_{i=1}^{r}\left\lceil{\frac{\sqrt{r+\rho_1s}(z_u^i-z_l^i)}{2\sqrt{\epsilon}}}\right\rceil N\right)$, where $N$ is the complexity to solve Problem (\ref{ODP:Bounded-CR}).

\section{Numerical experiments}\label{sec:Experiment}

In this section, we present  computational results of both the SCO algorithm (Algorithm~\ref{SCO}) and the SCOBB algorithm (Algorithm~\ref{SCOBBA}) for the OPD Problem (\ref{OPDP}). The algorithms are coded in Matlab  R2013b and run on a PC (3GHz, 16GB RAM).  All the linear and convex quadratic subproblems
in the algorithms are solved by the QP solver in Gurobi Optimizer 9.0.2  with Matlab interface (cf. \cite{gurobi}).

We apply our SCO and SCOBB algorithms to the OPD Problem (\ref{OPDP}) with  both synthetic data and historical stock data from NASDAQ\footnote{All the data  used in Section~\ref{sec:Experiment} can be downloaded on https://github.com/hezhiluo/OPD.}.
We compare our algorithms with the benchmark commercial software BARON (cf. \cite{int:Sahinidis96}). BARON, the acronym of Branch-And-Reduce Optimization Navigator, is a global optimization package for non-convex optimization problems.
To measure the effectiveness of the SCO algorithm, the relative  gap is computed by using the following formula
\begin{eqnarray*}\label{gap}
{\rm gap}: = (Opt.val -Obj.val)/\max\{1,|Opt.val|\}, %\label{eq-gap}
\end{eqnarray*}
where ``$Opt.val$'' denotes the global optimal value, ``$Obj.val$'' denotes the objective value at the solution found by the SCO algorithm.

In our numerical experiments, the stopping parameter $\epsilon$ is set at   $\epsilon=10^{-5}$. The maximum computational time is set as 3600 seconds. We document the notations in our computational results in Table 1.
\begin{table}[!h]\vspace*{-0.2in} \centering   \small
\caption{Notations in our numerical experiments.} \label{notation}
\begin{tabular}{l|l} \hline
$m$ & the number of variables\\
$s,q,r$&the rank of matrices $B^-$, $A^-$ and $B^-+A^-$, respectively, where $B^-$ and $A^-$ are given in (\ref{B-decomp})\\
Opt.val &the average global optimal value obtained by the algorithm    for 10 test instances\\
Time &the average CPU time of the algorithm   (unit: seconds)  for 10 test instances\\
Iter& the average number of iterations in the main loop of the algorithm   for 10 test instances\\
%Obj.val& the average objective value at the solution found by SCO for 10 test instances \\
$N_{\rm SCO}$& the average number of restarting SCO in the main loop of SCOBB for 10 test instances \\
$T_{\rm SCO}$ & the average CPU time of SCO used in SCOBB (unit: seconds) for 10 test instances\\
Gap& the average relative gap obtained by SCO for 10 test instances\\
\hline
\end{tabular}
\end{table}

\subsection{Numerical experiment with real data}\label{sect:tprd}
 In this subsection, we conduct numerical experiments using historical data from NASDAQ to illustrate the application of our algorithm. As discussed before, the most important parameters in an optimal portfolio deleveraging problem are the permanent and temporary price impact matrices. We thus first apply a linear regression model on historical stock data from NASDAQ to estimate the permanent and temporary price impact matrices, and then apply our algorithm to solve the optimal portfolio deleveraging problem on these stocks.

 We first describe the data set used in this numerical experiment. We use NASDAQ TotalView-ITCH data, which contain message data of order events.\footnote{Data are provided by LOBSTER (https://lobsterdata.com/).} The database documents all the order activities resulting in an update of the volumes at the best ask and bid prices, and these activities include visible orders' submissions, cancellations and executions. For each activity, the data set documents the price, trade sign (buy or sell), volume and order identity number. The timestamp of these activities is measured in seconds with decimal precision of at least milliseconds. We consider the optimal deleveraging problem on a portfolio consisting of six representative stocks from S\&P 100 components and five different sectors: Apple Inc.(AAPL), General Electric Company (GE), General Motors Company (GM), Coca-Cola Company (KO), PepsiCo, Inc. (PEP) and Walmart Inc. (WMT).

 \begin{table}[!ht]\caption{Descriptive statistics of the six stocks in our portfolio on Aug. 1st, 2018.}\label{tab:statistics}
 \centering
 \begin{tabular}{c|c|c|c}
 \hline
 \multirow{2}{*}{stock symbol} & aver. bid price & aver. hourly trade vol. &aver. trade size\\
 &(in dollar) &(in million) &(in share)\\
  \hline
 GM&\$37.2036 &0.1676 &111.2406\\
 \hline

 KO&\$46.3190&0.1157&124.9356\\
 \hline
 PEP&\$113.5936&0.1251&71.7464\\
 \hline
 WMT&\$88.3519&0.1539&89.7876\\
 \hline
 AAPL&\$200.2391&1.4356&108.8619\\
 \hline
  GE&\$13.3137&0.3808&645.9435\\
 \hline
 \end{tabular}
 \end{table}

% \begin{table}[!ht]\caption{Descriptive statistics of the six stocks in our portfolio on Aug. 1st, 2018.}\label{tab:statistics}
% \centering
% \begin{tabular}{c|c|c|c|c|c|c}
% \hline
%stock symbol &GM& KO & PEP &WMT&AAPL&GE\\
%  \hline
%\multirow{2}{*}{} aver. bid price&\$37.2036 &0.1676 &111.2406&\$37.2036 &0.1676 &111.2406\\
% (in dollar)&&&&&\\
% \hline
%%aver. hourly trade vol. &\$46.3190&0.1157&124.9356\\
%%  &(in million)
%% \hline
%%aver. trade size &\$113.5936&0.1251&71.7464\\
%%&(in share)
%% \hline
% \hline
% \end{tabular}
% \end{table}
 In particular, we first estimate the temporary and permanent price impact matrices of this portfolio with the data from 10:00 a.m. to 4:00 p.m. on Aug. 1st, 2018\footnote{We also consider different dates and up to twenty different stocks. The results are similar.}. We then apply our algorithm to solve the optimal portfolio deleveraging problem with the estimated matrices. See Table \ref{tab:statistics} for basic descriptive statistics of these six stocks on Aug. 1st, 2018.

 \textbf{Multiple linear regression model.} We apply the following multiple linear regression model to estimate the price impact matrices:
$$p_{it} = p_{i0} + c_{i} + \sum_{j=1}^{n}\gamma_{ji}(x_{jt}-x_{j0}) +\sum_{j=1}^{n} \lambda_{ji}y_{jt},~i=1,\cdots,6,$$
 where $(x_{jt}-x_{j0})$ is the cumulative signed trade volumes of the $j$th stock during the time period $[0,t)$ and $y_{jt}$ is the signed trade volume of the $j$th stock during a small time period $[t-1, t)$, $t\leq T$. The constant $c_{i}$ and impact parameters $\gamma_{ji}$ and $\lambda_{ji}$ are estimated by the linear regression model. Note that although $y_{jt}$ is the trading rate of a continuous trading in our model, it is impossible to trade continuously in practice. We thus consider the trade volume during a relatively small time period as $y_{jt}$ in our experiment. The small time period is set to be ten seconds in this estimation, while the whole time horizon $T$ is set to be 20 minutes. In other words, we take the impact of trades in last ten seconds on the price as an estimation of the temporary price impact, while the impact of trades in up to last twenty minutes on the price as an estimation of the permanent one.

 As a result, in our regression model, each data point is corresponding to one small time period, which is 10 seconds. With the data from 10:00 a.m. to 4:00 p.m. on Aug. 1st, 2018, there are $2160$ data points for each linear regression, and $18$ time horizons, each of which corresponds to 20 minutes. More specifically, each data point includes 1) the difference between the bid price vector at the end of the small time period and that at the beginning of the time horizon containing this small time period, 2) the cumulative signed trade volume vector during the time horizon, and 3) the signed trade volume vector during the small time period. We then document the temporary and permanent price impact matrices in Table \ref{tab-IM}, where  all the temporary and permanent impact coefficients have been multiplied by $10^4$.  The statistics of this linear regression are presented in the e-companion (Table~\ref{tab:estimationstat}). In particular, for each of the six linear regressions, we document the R-square and p-value of the F-test, which show a significant linear regression relationship between the response and the predictor variables.

\begin{table}[!ht]
\caption{The price impact matrices estimated from our sample data on Aug. 1st, 2018. The $ij$th element of the matrix denotes the temporary price impact parameter of the asset $i$ on the asset $j$, multiplied by $10^4$.}
\label{tab-IM}
\begin{center}
\begin{tabular}{c|cccccc}
\hline
\multicolumn{7}{c}{The temporary price impact matrix.}\\
\hline
 &GM &KO &PEP&WMT&AAPL& GE\\
\hline
GM&0.081784     & -0.0039955   & -0.0082433 & -0.027238  & 0.018445    &  0.0017027  \\
KO&-0.0076704   & 0.049933    &   0.079503   & 0.058901     & -0.020779     &  -0.0010544  \\
PEP&  -0.0092156  & 0.028484  &  0.074409    &  0.011151   &   -0.052514   & -0.0071483     \\
WMT&  -0.0085645 & 0.013953   &  -0.006123    &  0.10158     &  0.1043   &  -0.008582  \\
AAPL&  0.00043804 & -5.49E-05 &  2.41E-05    &   0.0017344      &  0.017789     & 0.00017018   \\
GE& 0.0045241  &0.00086434   & 0.0021353 &   0.0095049   & -0.0049313    &   0.0032635 \\
\hline
\multicolumn{7}{c}{The permanent price impact matrix.}\\
\hline
&GM &KO &PEP&WMT&AAPL& GE\\
\hline
GM&0.10236    & -0.0061319     & -0.019132  & -0.0073272   & 0.22656    & -2.42E-05 \\
KO&-0.004774 &  0.039061       &  0.091118   &0.050312    &  0.058154  &  -0.00145 \\
PEP&-0.0084035  & 0.059625     & 0.085173    & -0.0082462  &  -0.15316   &  -0.0067892   \\
WMT&0.008701 & 0.0084812       & 0.0034464   &  0.066367  &  0.041291   &  -0.00388  \\
AAPL& -0.0018108 &-0.00051277  &  -0.001625  &  -0.0044423  &  0.014685     & -0.00020722  \\
GE&  0.0060147  & -0.0024202   &3.46E-05     & 0.0062422   &  -0.0052475  &  0.0039042  \\
\hline
\end{tabular}
\end{center}
\end{table}
%\begin{table}[!ht]
%\caption{The permanent price impact matrix estimated from our sample data on Aug. 1st, 2018. The $ij$th element of the matrix denotes the permanent price impact parameter of the asset $i$ on the asset $j$, multiplied by $10^4$.}
%\label{tab-PPIM}
%\begin{center}
%\begin{tabular}{c|cccccc}
%\hline
% &GM &KO &PEP&WMT&AAPL& GE\\
%\hline
%GM&0.10236    & -0.0061319     & -0.019132  & -0.0073272   & 0.22656    & -2.42E-05 \\
%KO&-0.004774 &  0.039061       &  0.091118   &0.050312    &  0.058154  &  -0.00145 \\
%PEP&-0.0084035  & 0.059625     & 0.085173    & -0.0082462  &  -0.15316   &  -0.0067892   \\
%WMT&0.008701 & 0.0084812       & 0.0034464   &  0.066367  &  0.041291   &  -0.00388  \\
%AAPL& -0.0018108 &-0.00051277  &  -0.001625  &  -0.0044423  &  0.014685     & -0.00020722  \\
%GE&  0.0060147  & -0.0024202   &3.46E-05     & 0.0062422   &  -0.0052475  &  0.0039042  \\
%\hline
%\end{tabular}
%\end{center}
%\end{table}

As shown in Table \ref{tab-IM}, the cross impacts between stocks in the same sector (PEP and KO) are positive, while those cross-sector impacts could be negative. This finding is consistent with that by \cite{pasquariello2013strategic}. In addition, the impact of the trades of one stock on the price of one another stock, although much smaller than the impact of the trades on its own price, is non-ignorable.  This finding is consistent with our intuition: due to possible portfolio rebalancing trades, strategical trading activities of sophisticated speculators and some correlated liquidity shocks, the prices and trades of different stocks may have co-movements. The fact that these estimated temporary and permanent price impact matrices are not diagonal also demonstrates the meaningfulness of our algorithm considering non-diagonal price impact matrices.

Moreover, the impact matrices are not only non-diagonal, but also asymmetric, which is also found in the empirical results documented in Table V by \cite{pasquariello2013strategic}. The intuition is that the relationship between these fundamentally uncorrelated stocks could be complicated and that stocks may vary in the liquidity level: a more liquid stock could be less significantly affected by trades.\footnote{We want to point out that, although in the theoretical models without the consideration of bid-ask spread, the price impact matrices are usually assumed to be symmetric and positive semi-definite to avoid arbitrage possibilities, an asymmetric price impact matrix with limited number of negative eigenvalues does not necessarily bring in arbitrage opportunities in the real market with a positive bid-ask spread.}
%\begin{table}[!ht]\caption{The statistics of the linear regression model for each stock.}\label{tab:estimationstat}
%\centering
%\begin{tabular}{c|c|c|c|c|c|c}
%\hline
%stock symbol &GM &KO &PEP&WMT&AAPL& GE\\
%\hline
%p-value&0.000&0.000&0.000&0.000&0.000&0.000\\
%\hline
%R-square&51.10\% &63.89\% &56.41\% &40.09\% &47.78\% &48.09\%\\
%\hline
%\end{tabular}
%\end{table}

\begin{example}\label{exp_nasdaq} We consider the optimal deleveraging Problem (\ref{OPDP}) with the price impact matrices $\Lambda$ and $\Gamma$ estimated and shown in Table \ref{tab-IM}.  Based on the descriptive statistics in Table \ref{tab:statistics}, we set the initial holding $x_{0}$ and the initial price $p_{0}$ in Problem \eqref{OPDP} as
\[x_{0}=(2000,2000,2000,2000,8600,5000)^{T},\quad p_{0}=(37.39,46.42,113.91,88.94,198.76,13.4)^{T}.\]
\end{example}
Suppose that the initial liability is $l_0=2262000$ and the required leverage ratio is $\rho_1=18$. One can easily calculate that the initial equity is $e_0=87656$, and the initial liability equity
ratio  is $l_0/e_0=25.8054$, and testify that Assumption~\ref{asmp1} holds.

We consider the relaxed problem
by removing the leverage constraint in Problem \eqref{OPDP}:
\begin{eqnarray}\label{ODP-nlc}
 \max_{y\in{\mathbb R}^m}\left\{e_1(y):-x_0\le y\le 0\right\}.
 \end{eqnarray}
We solve Problem~(\ref{ODP-nlc}) by the ADMBB algorithm in \cite{lblp19} to identify the global optimal solution $\hat{y}$.
Denote $\rho_{\max}:=l_1(\hat y)/e_1(\hat y)$. The required leverage ratio should be less than $\rho_{\max}$ to be ``active'' in the problem.
%In addition, without the leverage constraint \eqref{eqn:leverage-constraint2}, the relaxed version of optimization problem \eqref{OPDP} would
%yield an optimal solution with the corresponding leverage ratio at time $t=1$ being $\rho_{max} = 25.42 > \rho_{1}$.
Solving Problem~(\ref{ODP-nlc}) associated with the above data yields $\rho_{\max} = 25.42 > \rho_{1}$.
In other words, our problem setting guarantees that the trader in this example is capable to meet his liability obligation by liquidating and is forced to liquidate by the leverage constraint.

According to the discussion in Section~\ref{sec:formulation-SD}, we first reformulate the OPD problem as Problem (\ref{OPDP-SD}) with $s=1, r=3$. We then apply our SCO and SCOBB algorithms and BARON\footnote{In this example, BARON cannot solve Problem \eqref{OPDP-sym} within its default time limit, 500 seconds. We thus apply BARON to Problem \eqref{OPDP-SD}, instead of Problem \eqref{OPDP-sym} in this example.} to Problem \eqref{OPDP-SD}. From Proposition~\ref{prop-SD}, the optimal solution to the original problem  can be obtained from the solution to Problem \eqref{OPDP-SD}.

%Then, we apply our SCO and SCOBB algorithms and BARON to solve problem (\ref{OPDP-SD}) associated with the above data.
Numerical results of SCO, SCOBB and BARON are reported in Table~\ref{tab_exam},  where ``Opt.solution'' denotes  the  optimal solution $y^*$ to Problem (\ref{OPDP}) derived from the solution by the algorithms.  In particular, as shown in Table~\ref{tab_exam}, our SCOBB algorithm  can find  the global optimal solution within 0.1 second.
%The globally optimal solution  is successfully identified by the SCOBB algorithm with 31 iterations with CPU time 0.38 second.
%Thus, the globally optimal solution of problem~(\ref{OPDP}) found by the SCOBB is
%The corresponding optimal trading rates  are
%\[y^* =(-1478.8137, -446.7456, 0, 0, -2754.7015, -5000)^T \]
%The optimal value, which is the equity after trading, is thus $e_1(y^*)=87523.223953$.
The optimal value identified by our SCOBB algorithm is the same as that by BARON. In addition, our SCO algorithm, which is proven to converge to a KKT point, also finds the global optimal solution.
%provides a good approximation for the global optimizer.
The result also shows that the leverage ratio constraint $g(y^*)\leq 0$ could be active, even in this case where the price impact matrices are asymmetric with some negative elements.
\begin{table}[!ht] \vspace*{0.0in}
\caption{\label{tab_exam} Numerical results of our algorithms and BARON for Problem (\ref{OPDP}) in Example~\ref{exp_nasdaq}.}
\begin{center}
\small
\begin{tabular}{c|c|c|c|c|c } \hline
 Algorithm  & Opt.val  & Time &  Iter & Opt.solution   & $g(y^*)$   \\
\hline
 SCO&87523.223953&0.03&8&$(-1478.8137, -446.7456, 0, 0, -2754.7015, -5000)^T$&-0.000520\\
  \hline
SCOBB &87523.223953 &0.07 &4&$(-1478.8137, -446.7456, 0, 0, -2754.7015, -5000)^T$ &-0.000520 \\
\hline
BARON&87523.223953 & 0.37 &91&$(-1478.8137, -446.7456, 0, 0, -2754.7015, -5000)^T$ &0.000000\\
 \hline
 \end{tabular}
 \end{center}
 \vspace{2mm}
\end{table}

To illustrate the effects of cross impacts on selling priority discussed in Proposition~\ref{prop-priority}, we furthermore consider an example where the initial price $p_{0}$ is set as $p_{0} = (40,40,100,100,200,20)^{T}.$ The price impact matrices, initial holding $x_{0}$ and initial liability $l_{0}$ are the same as Example~\ref{exp_nasdaq}. Similarly, one can easily calculate that the initial leverage ratio in such an example is $l_0/e_0 = 19.17$, the leverage ratio yielding by Problem~(\ref{ODP-nlc}) is $\rho_{max}=18.75$ and that Assumption~\ref{asmp1} holds.
%Solving Problem~(\ref{ODP-nlc}) associated with these data yields $\rho_{\max}=18.75$.
We thus set the required leverage ratio $\rho_{1}$ varying from 16 to 8, to show the selling priority in this example and to make sure that the trader is forced to liquidate and is capable to meet the leverage requirement. The numerical results of SCOBB are presented in Table~\ref{tab2_exam3}.
\begin{table}[!ht] \vspace*{0.0in}
\caption{Numerical results of SCOBB for Problem (\ref{OPDP}) in Example~\ref{exp_nasdaq} with $p_0=(40,40, 100,100,200,20)^T$.}
\label{tab2_exam3}
\begin{center}
\begin{tabular}{c|c|c|c|c|c } \hline
 $\rho_1$ & Opt.val & Time & Iter & Opt.solution $y^*$ & $g(y^*)$ \\
  \hline
 8 & 117785.128650 & 0.03 & 1 &$(-2000, -2000, 0, 0, -5313.1975, -5000 )^T$& -5.2e-5 \\
 10 & 117815.903048 & 0.11 & 2 &$(-2000, -1362.9708, 0, 0, -4263.8450, -5000)^T$& -3.0e-5 \\
 12 & 117848.654786 & 0.04 & 2 &$(-1999.9989, -581.5282, 0, 0, -3241.9450, -4999.9992)^T$& -3.5e-3 \\
 14 & 117887.549258 & 0.08 & 5 &$(-1053.5473, -56.2028, 0, 0, -2344.8658, -5000)^T$& -1.5e-4 \\
 16 & 117938.299432 & 0.06 & 6 &$(-0.0004, 0, 0, 0, -1371.3264, -5000)^T$& -4.8e-5  \\
 \hline
\end{tabular}
\end{center}
\end{table}
One can observe that the first asset is prioritized for selling, compared with the second asset, while the initial price and holding of these two assets are the same. Although the four conditions in Proposition~\ref{prop-priority}, which provide a sufficient, but not necessary condition for prioritized sales, do not hold in this case. Inspired by conditions (ii) and (iv), we can easily testify that the first asset (GM) is less positively correlated with the other assets than the second one (KO) by
\begin{eqnarray*}
&\sum\limits_{\substack{k\neq 1,2\\ k=1}}^{6}\hat{\Lambda}_{1k}-\sum\limits_{\substack{k\neq 1,2\\ k=1}}^{6}\hat{\Lambda}_{2k} = -0.0294<0,~ \sum\limits_{\substack{k\neq 1,2\\ k=1}}^{6}\hat{\Gamma}_{1k}-\sum\limits_{\substack{k\neq 1,2\\ k=1}}^{6}\hat{\Gamma}_{2k} = -0.0940<0,&\\
&\sum\limits_{k=1}^{6} (\gamma_{k1} - \gamma_{1k}) = -0.1942 < \sum\limits_{k=1}^{6}(\gamma_{k2} - \gamma_{2k})=-0.1343.&
\end{eqnarray*}
As a result, this illustrative example demonstrates that the cross impacts are important in the deleveraging problem and the selling priority inside the portfolio.

\begin{example}\label{exam2_nasdaq}  We consider the OPD problems with portfolios consisting  of up to twenty assets from NASDAQ during one sample month and with the historical data from the same database as Example 3. The describe statistics of the assets are documented in our e-companion (See Table~\ref{tab:example4-statistics}).  We estimate the impact parameters in a similar way as described in Example~\ref{exp_nasdaq}. The price impact matrices and related statistics are available upon request. In particular, we consider ten different sets of data, and  obtain ten instances of OPD problem for each portfolio size $m=10,15,20$.
\end{example}

We assume that the initial liability equity ratio is $l_0/e_0=25$ and  the required leverage ratio is $\rho_1=18$, and testify that Assumption~\ref{asmp1} holds for all instances. As in Example~\ref{exp_nasdaq}, we transform the OPD problem into Problem~(\ref{OPDP-sym}) and then Problem~(\ref{OPDP-SD}). We then apply our SCOBB and SCO algorithms separately to Problem~(\ref{OPDP-SD}), and BARON to Problem~(\ref{OPDP-sym}). The optimal solution to the original problem can be  obtained from the solution to the transformed problems by using Proposition~\ref{prop-SD}. According to the identified optimal solutions, the leverage ratio constraint $g(y^*)\leq 0$ could be inactive for some instances. That is to say, traders might sell more than the leverage requirement in the OPD problem, leading to $g(y^*)< 0$. As we discussed in Section~\ref{sec:property}, the intuition is that there might be some hedge properties inside the portfolio, such that the trader could utilize non-zero cross impact among different assets.

%\begin{table}[!ht] \vspace*{0.0in}
%\caption{Average rank of small-scaled instances of OPD problem with historical data in Example~\ref{exam2_nasdaq}, and the average performances of our SCOBB and SCO algorithms and BARON. }
%\label{tab1-opdp-real-20}\footnotesize
%\begin{center}
%\begin{tabular}{|c|ccc|crrcc|cr|ccc|} \hline
%\multicolumn{1}{|c|}{Size}& \multicolumn{3}{c|}{Rank}   & \multicolumn{5}{c|}{SCOBB}& \multicolumn{2}{c|}{BARON} & \multicolumn{3}{c|}{SCO}    \\
%\hline
%$m$ & $s$& $q$ &$r$ & Opt.val & Time & Iter  & $N_{sco}$& $T_{sco}$& Opt.val & Time  & Gap & Time& Iter     \\
%\hline
%%5 & 0.40 & 0.1 & 0.50 &   &   &      &  &  &  &   &   &    \\
%10 & 1.7 & 0.8 & 2.5 & 1701524.8172 & 0.57 & 34.9 & 0.9 & 0.15    & 1578851.10(7) & 0.55   & 2.4e-7 & 0.12 & 8.7    \\
%15 & 3.1 & 2.0 & 5.1 & 2859387.0511 & 2.70 & 191.9 & 0.2 & 0.23   & 3074218.44(7) & 13.90  & 1.0e-8 & 0.18 & 15.3   \\
%20 & 4.8 & 3.4 & 8.2 & 3799765.8138 & 20.27 & 1299.2 & 0.1 & 0.38 & 3889242.28(8) & 112.25 & 1.0e-8 & 0.37 & 25.5   \\
%\hline
%\end{tabular}
% \end{center}
% \vspace{2mm}
%{\bf Remark:} The number in parentheses in the column of ``Opt.val'' for BARON stands for the number of the instances for which BARON  can verify the global optimality of the  solution within 500 seconds. ``Time" and ``Opt.val"  for BARON denote  the average CPU time and optimal value  for the instances  that are globally solved by BARON in the ten instances, respectively.
%\end{table}

\begin{table}[!ht] \vspace*{0.0in}
\caption{Average rank of small-scaled instances of OPD problem with price impact matrices estimated from the historical data in Example~\ref{exam2_nasdaq}, and the average performances of our SCOBB and SCO algorithms and BARON. }
\label{tab1-opdp-real-20}\footnotesize
\begin{center}
\begin{tabular}{|c|ccc|c|crrcc|cr|ccc|} \hline
\multicolumn{1}{|c|}{Size}& \multicolumn{3}{c|}{Rank} & \multicolumn{1}{c|}{}    & \multicolumn{5}{c|}{SCOBB}& \multicolumn{2}{c|}{BARON} & \multicolumn{3}{c|}{SCO}    \\
\hline
$m$ & $s$& $q$ &$r$ &$\rho_{\max}$& Opt.val & Time & Iter  & $N_{sco}$& $T_{sco}$& Opt.val & Time  & Gap & Time& Iter     \\
\hline
10 & 1.7 & 0.8 & 2.5 & 21.90& 1701524.8172 & 0.57 & 34.9 & 0.9 & 0.15    & 1578851.10(7) & 0.55   & 2.4e-7 & 0.12 & 8.7    \\
15 & 3.1 & 2.0 & 5.1 & 19.80& 2859387.0511 & 2.70 & 191.9 & 0.2 & 0.23   & 3074218.44(7) & 13.90  & 1.0e-8 & 0.18 & 15.3   \\
20 & 4.8 & 3.4 & 8.2 & 20.54& 3799765.8138 & 20.27 & 1299.2 & 0.1 & 0.38 & 3889242.28(8) & 112.25 & 1.0e-8 & 0.37 & 25.5   \\
\hline
\end{tabular}
 \end{center}
 \vspace{2mm}
{\bf Remark:} The number in parentheses in the column of ``Opt.val'' for BARON stands for the number of the instances for which BARON  can verify the global optimality of the  solution within 500 seconds. ``Time" and ``Opt.val"  for BARON denote  the average CPU time and optimal value  for the instances  that are globally solved by BARON in the ten instances, respectively.
\end{table}

 In Table~\ref{tab1-opdp-real-20}, we summarize the average performances of SCOBB, SCO and BARON on these small-scaled instances of OPD problem with the historical data on these twenty assets from NASDAQ.
 As shown in Table \ref{tab1-opdp-real-20}, our SCOBB algorithm can identify the global optimal solution for all test instances within around 20 seconds, while BARON can only solve 22 out of 30 instances within its default computational time (500 seconds). Moreover, it takes much longer time for BARON than SCOBB to identify the solution to these solved instances with $m=15,20$. We also observe that  the KKT point found  by the SCO algorithm  provides a good approximation for the global optimizer in terms of the gap.

In addition, the computational time of our SCOBB algorithm grows exponentially in terms of the number of branches in the B\&B framework, and thus the number of negative eigenvalues of the matrices $\hat{\Lambda}-\frac{1}{2}\hat{\Gamma}$ and $\hat{\Lambda}+\frac{1}{2}\hat{\Gamma}$ in the objective function and the constraint. As a result, we document the average numbers of negative eigenvalues in these instances with historical data in Table~\ref{tab1-opdp-real-20}. The numerical results demonstrate that the total number of negative eigenvalues of the two matrices, $r$ defined in \eqref{SD}, is usually not too large in the real problem. This fact guarantees that our SCOBB algorithm can solve OPD problem within short computational time, and thus serve a good purpose as a practical algorithm for the investors.

\subsection{Numerical results for randomly generated test problems}\label{sect:rgtp}
In this subsection, we test the SCOBB algorithm on the OPD Problem (\ref{OPDP}) with randomly
generated instances. %The computational time of our SCOBB algorithm grows exponentially in terms of the number of branches in the B\&B framework, and thus the number of negative eigenvalues of the matrices $\hat{\Lambda}-\frac{1}{2}\hat{\Gamma}$ and $\hat{\Lambda}+\frac{1}{2}\hat{\Gamma}$ in the objective function and the constraint.
As shown in Example~\ref{exam2_nasdaq}, the number of negative eigenvalues is limited in the real problem. Our experiments are thus restricted to instances with a few negative eigenvalues, i.e., a small $r$ defined in \eqref{SD}.

Recall that \cite{Brown10} assume that the matrix $\Lambda-\frac{1}{2}\Gamma$ is positive definite and that both the permanent and temporary matrices are diagonal and positive definite. Under their assumption, the number of negative eigenvalues is $r=0$. As stated in \cite{Brown10}, the restriction on the matrix $\Lambda-\frac{1}{2}\Gamma$ actually prevents the trader from trading infinite size to obtain arbitrarily large equity, and thus ensures that the deleveraging problem is well posed. Our global algorithm SCOBB is also customized for the problem with small $r$.
 %We will show in the following numerical results that the performances of our SCOBB algorithm are still much better than the BARON package for the small-scale problem with $r$ up to 15.

\begin{example}\label{exp_limitr}
Inspired by Example~\ref{exp_nasdaq} with real data, we randomly generate the parameters $(\Lambda,\Gamma,p_0,x_0)$ in our test instances as follows. The entries of $p_0$ are drawn from $U[10,100]$, i.e., uniformly distributed within interval $[10,100]$, and entries of $x_0$ are drawn from $U[500,1000]$. The matrices $\hat{\Lambda}-\frac{1}{2}\hat{\Gamma}$ and $\hat{\Lambda}+\frac{1}{2}\hat{\Gamma}$ are randomly generated\footnote{The way we generate matrices with an arbitrary fixed number of negative eigenvalues is documented in our e-companion.} to have exactly $s$ and $q$ negative eigenvalues, respectively, with elements lying inside $[10^{-6},10^{-5}]$.
% Thus, $\hat{\Lambda}-\frac{1}{2}\hat{\Gamma}$ and $\hat{\Lambda}+\frac{1}{2}\hat{\Gamma}$ have exactly $s$ and $q$ negative eigenvalues, respectively.
We consider the optimal deleveraging problem with an initial liability equity ratio  $l_0/e_0=25$ and  the required leverage ratio $\rho_1=18$. Assumption~\ref{asmp1} is satisfied for each random test instance.
We then apply our SCO and SCOBB algorithms and BARON to randomly generated ten test instances for varied fixed sizes and numbers of negative eigenvalues.

\begin{table}[!ht] \vspace*{0.0in}
\caption{Comparison among the average performances of our SCOBB and SCO algorithms and BARON for small-scaled instances of OPD problem with $\rho_1=18$. }
\label{tab1-opdp}\footnotesize
\begin{center}
\begin{tabular}{|cccc|c|crrc|cc|cccc| } \hline
\multicolumn{4}{|c|}{Size}  & \multicolumn{1}{c|}{} & \multicolumn{4}{c|}{SCOBB}& \multicolumn{2}{c|}{BARON} & \multicolumn{4}{c|}{SCO}    \\
\hline
$m$ & $s$& $q$ &$r$ &$\rho_{\max}$& Opt.val & Time & Iter  & $\frac{l_1(y^*)}{e_1(y^*)}$ & Opt.val & Time  & Gap & Time& Iter   &  $\frac{l_1(y^*)}{e_1(y^*)}$   \\
\hline
20 & 2 & 3 & 5  & 23.10 & 32101.20113 & 0.37 & 24.8    & 18   & 29804.77150(1) & 180.09& 0 & 0.08 & 14.5  &  18    \\
20 & 4 & 4 & 8  & 22.93 & 31192.74541 & 1.53 & 129.4   & 18   & 27422.59496(2) & 245.43& 0 & 0.08 & 16.9  &  18  \\
20 & 5 & 5 & 10 & 22.93 & 31329.06380 & 2.83 & 235.8   & 18   & 30025.98130(6) & 397.78& 0 & 0.11 & 21.8  &  18  \\
20 & 7 & 8 & 15 & 22.58 & 33949.10610 & 47.83 & 3063.5 & 18   & 33233.46897(4) & 420.41& 0 & 0.11 & 22.0  &  18  \\
 \hline
30 & 2 & 3 & 5  & 23.46 & 47919.2267 & 0.40 & 19.1     & 18   & $-$ & $-$ & 0 &0.09 & 12.8 &  18   \\
30 & 4 & 4 & 8  & 23.53 & 50733.4138 & 3.84 & 193.3    & 18   & $-$ & $-$ & 0 &0.17 & 19.4 &  18  \\
30 & 5 & 5 & 10 & 23.34 & 46827.9257 & 4.62 & 243.6    & 18   & $-$ & $-$ & 0 &0.13 & 16.4 &  18    \\
30 & 7 & 8 & 15 & 23.26 & 48078.3605 & 53.74 & 2498.4  & 18   & $-$ & $-$ & 0 &0.16 & 22.3 &  18   \\
\hline
 \end{tabular}
 \end{center}
  \vspace{2mm}
{\bf Remark:} The number in parentheses in the column of ``Opt.val'' for BARON stands for the number of the instances for which BARON  can verify the global optimality of the  solution within 1000 seconds. ``Time" and ``Opt.val"  for BARON denote  the average CPU time and optimal value  for the instances  that are globally solved by BARON   in the ten instances, respectively. The sign ``$-$" stands for the situations where the method fails to identify the global solution within 1000 seconds in all the ten instances.
\end{table}

We first compare the performances of SCOBB, SCO and  BARON for small-scale instances of Problem (\ref{OPDP}), i.e., a small-scale portfolio with the number of assets being $m=20$ or $m=30$, with $r=5,8,10,15$. The results documented in Table~\ref{tab1-opdp} show that the SCOBB algorithm is able to find the global optimal solution within 100 seconds for all the test instances, while  BARON  can only identify the global optimal solution in 13 out of 80 test instances within 1000 seconds. Also, BARON requires more CPU time than SCOBB for the solved instances. In addition, for most of test instances, BARON only reported the best solution obtained within 1000 seconds and failed to verify the global optimality of the obtained solution. The computational time of BARON grows rapidly in terms of the number of assets. As shown in Table~\ref{tab1-opdp}, BARON fails to identify the global optimal solution for all the instances with $m=30$.
Meanwhile, the CPU time of our SCOBB algorithm grows rapidly in terms of $r$, which results from the branch-and-bound framework.
%\chen{Note that the existence of negative eigenvalues, i.e., $r>0$, would bring in multiple local solution traps for the local optimization algorithm. (Yuanyuan: please double check this sentence in case I explain it wrong.) }
The result shows that our SCO algorithm always provides the global optimal solutions for all the instances in our numerical experiments.

 \begin{table}[!ht] \vspace*{0.0in}
\caption{Average performances of SCOBB and SCO for medium- and large-scale OPD instances  with $\rho_1=18$.}
\label{tab2-opdp}  \centering\small
\begin{tabular}{|cccc|c|crrc|crcc| } \hline
\multicolumn{4}{|c|}{Size}& \multicolumn{1}{c|}{}   & \multicolumn{4}{c|}{SCOBB} &\multicolumn{4}{c|}{SCO}    \\
\hline
$m$& $s$ &$q$ &$r$   &$\rho_{\max}$&Opt.val   & Time  & Iter  & $\frac{l_1(y^*)}{e_1(y^*)}$   &Gap  & Time & Iter & $\frac{l_1(y^*)}{e_1(y^*)}$ \\
\hline
100 & 3 & 3 & 6 & 24.16 & 154874.9286 & 4.66 & 43.5     &  18 &   0 & 0.98 & 18.7 & 18 \\
100 & 4 & 4 & 8 & 23.93 & 163126.9455 & 13.12 & 135.5   &  18 &   0 & 1.21 & 20.9 & 18  \\
100 & 5 & 5 & 10& 23.98 & 159477.5510 & 31.16 & 351.9   &  18 &   0 & 0.99 & 20.1 & 18 \\
\hline
200 & 3 & 3 & 6 & 24.20 & 318277.9335 & 28.83 & 52.5    & 18 &   0  & 7.60 & 25.7 &18   \\
200 & 4 & 4 & 8 & 24.23 & 311463.5154 & 47.39 & 112.6   & 18 &   0  & 4.52 & 14.8 &18  \\
200 & 5 & 5 & 10& 24.20 & 308724.1236 & 131.46 & 324.7  & 18 &   0  & 6.05 & 21.1 &18  \\
\hline
300 & 3 & 3 & 6 & 24.35 & 468396.8383 & 79.63 & 60.4     & 18& 0 & 21.41 & 22.3  & 18 \\
300 & 4 & 4 & 8 & 24.33 & 467601.3706 & 159.67 & 139.7   & 18& 0 & 24.76 & 26.3  & 18 \\
300 & 5 & 5 & 10& 24.33 & 467005.5263 & 306.24 & 297.0   & 18& 0 & 24.33 & 25.8  & 18  \\
\hline
400 & 3 & 3 & 6 & 24.38 & 615973.2506 & 168.32 & 54.9   & 18  & 0  & 24.89 & 13.3& 18  \\
400 & 4 & 4 & 8 & 24.43 & 610797.3286 & 285.05 & 92.6   & 18  & 0  & 57.00 & 28.3& 18   \\
400 & 5 & 5 & 10& 24.41 & 622042.8137 & 999.60 &355.8   & 18  & 0  & 44.33 & 20.2& 18  \\
\hline
500 & 3 & 3 & 6 & 24.44 & 776512.5791 & 435.50 & 54.4    & 18& 0  & 97.18  & 22.5 &18\\
500 & 4 & 4 & 8 & 24.45 & 787282.3202 & 986.39 & 130.1   & 18& 0 & 126.84  & 30.6 &18 \\
500 & 5 & 5 & 10& 24.43 & 783519.0066 & 1562.69 & 214.1  & 18& 0 & 135.15  & 26.4 &18  \\
\hline
\end{tabular}
\end{table}

Table~\ref{tab2-opdp} summarizes the average numerical results of  SCOBB and SCO for  medium- and large-scale instances of Problem (\ref{OPDP}) with  $r=6,8,10$. As one can see from Table~\ref{tab2-opdp} that  SCOBB  can effectively find the global optimal solution  for all the test instances of Problem (\ref{OPDP})  within 1600 seconds, with the number of assets up to $m=500$.
%\chen{In fact, the computational complexity of our algorithm is mainly affected by the rank $r$, but not the number of assets $m$. As we discussed at the beginning of this section, as long as there are limited number of negative eigenvalues $r$, our SCOBB algorithm can identify the global optimal solution within a short time and a few iteration steps, even for a large size portfolio. }
In addition, the SCO algorithm can always obtain the global optimal solutions in less CPU time for all the medium- and large-scale instances in our experiments.
\end{example}

\section{Conclusions}\label{sec:Conclusion}
Market impact is crucial for traders, especially large traders, to deleverage a portfolio during a short time period. Cross impact among the assets are non-ignorable, due to financial constraints or the portfolio rebalancing trades from sophisticated speculators. In this paper, we investigate an optimal deleveraging problem with cross impact. The objective is to maximize equity while meeting a prescribed leveraged ratio requirement. We obtain some analytical insights regarding the optimal deleveraging strategy. Especially, we find that with the cross price impact, traders may sell more, instead of precisely satisfy the leverage requirement, if there are hedge properties inside the portfolio. In addition, among all the assets inside the portfolio, the trader should sell more of the asset which is more liquid and less correlated with the other assets. These analytical properties are verified by our empirical and numerical experiments.

Our optimization model is a non-convex quadratic program with non-homogeneous quadratic and box constraints, which is NP-hard. We reformulate the non-convex quadratic program as a D.C. program via  spectral decomposition and simultaneous diagonalization, and then propose two efficient algorithms for it: the SCO method
and the SCOBB method. We show that our SCO algorithm converges to a KKT point of the transformed problem, while the SCOBB algorithm can efficiently identify the global $\epsilon$-optimal solution. We establish the global convergence of the SCOBB algorithm and estimate its complexity.

We estimate the price impact matrices with the historical data from NASDAQ in our numerical experiments, to demonstrate the application of our algorithms. Randomly generated instances are also considered. According to our numerical experiments, our SCOBB algorithm can effectively identify the global optimal solution to medium- and large-scale instances of the optimal deleveraging problem with limited number of negative eigenvalues of the matrices in the quadratic terms in our optimization problem. A future research topic is to investigate whether we can develop  effective global algorithms for general optimal deleveraging problem without this restriction. Meanwhile, although we cannot prove the global optimality, our SCO algorithm always provides the global optimal solution in our numerical experiments for all the instances within short computational time, and thus also serve the purpose as an efficient algorithm for the optimal deleveraging strategy.

\theendnotes
\ACKNOWLEDGMENT{%
The research of H. Luo and H. Wu is supported by NSFC grants 11871433 and 11371324  and the Zhejiang Provincial NSFC grants  LZ21A010003 and LY18A010011. The research of Y. Chen is supported by NSFC grants 72001105. The research of D. Li is supported by Hong Kong Research Grants Council under grants 14213716 and 14202017.
}

%\bibliographystyle{informs2014}
%\bibliography{AssetLiqid}

\begin{thebibliography}{36}
\providecommand{\natexlab}[1]{#1}
\providecommand{\url}[1]{\texttt{#1}}
\providecommand{\urlprefix}{URL }

\bibitem[{Almgren \protect\BIBand{} Chriss(2000)}]{Almgren00}
Almgren R, Chriss N (2000) Optimal execution of portfolio transactions.
  \emph{J. Risk} 3(2):5--39.

\bibitem[{Andrade et~al.(2008)Andrade, Chang, \protect\BIBand{}
  Seasholes}]{Andrade08}
Andrade SC, Chang C, Seasholes MS (2008) Trading imbalances, predictable
  reversals, and cross-stock price pressure. \emph{J. Financial Econom.}
  88(2):406--423.

\bibitem[{Audet et~al.(2000)Audet, Hansen, Jaumard, \protect\BIBand{}
  Savard}]{Audet00}
Audet C, Hansen P, Jaumard B, Savard G (2000) A branch and cut algorithm for
  nonconvex quadratically constrained quadratic programming. \emph{Math.
  Program.} 87:131--152.

\bibitem[{Ben-{T}al \protect\BIBand{} den Hertog(2014)}]{Ben-Tal14}
Ben-{T}al A, den Hertog D (2014) Hidden conic quadratic representation of some
  nonconvex quadratic optimization problems. \emph{Math. Program., Ser. A}
  143:1--29.

\bibitem[{Ben-{T}al \protect\BIBand{} Teboulle(1995)}]{Ben-Tal95}
Ben-{T}al A, Teboulle M (1995) Hidden convexity in some nonconvex quadratically
  constrained quadratic programming. \emph{Math. Program.} 72(1):51--63.

\bibitem[{Benzaquen et~al.(2017)Benzaquen, Mastromatteo, Eisler,
  \protect\BIBand{} Bouchaud}]{Benzaquen17}
Benzaquen M, Mastromatteo I, Eisler Z, Bouchaud JP (2017) Dissecting
  cross-impact on stock markets: An empirical analysis. \emph{J. Stat. Mech.
  Theory E. {\rm Online at stacks.iop.org/JSTAT/2017/023406}}
  \urlprefix\url{http://dx.doi.org/10.1088/1742-5468/aa53f7}.

\bibitem[{Bertsimas et~al.(1999)Bertsimas, Hummel, \protect\BIBand{}
  Lo}]{Bertsimas99}
Bertsimas D, Hummel P, Lo AW (1999) Optimal control of execution costs for
  portfolios. \emph{Comput. Sci. Engrg.} 1(6):40--53.

\bibitem[{Brown et~al.(2010)Brown, Carlin, \protect\BIBand{} Lobo}]{Brown10}
Brown DB, Carlin B, Lobo MS (2010) Optimal portfolio liquidation with distress
  risk. \emph{Management Sci.} 56(11):1997--2014.

\bibitem[{Burer \protect\BIBand{} Vandenbussche(2008)}]{Burer08}
Burer S, Vandenbussche D (2008) A finite branch-and-bound algorithm for
  nonconvex quadratic programming via semidefinite relaxations. \emph{Math.
  Program.} 113(2):259--282.

\bibitem[{Burer \protect\BIBand{} Vandenbussche(2009)}]{Burer09}
Burer S, Vandenbussche D (2009) Globally solving box-constrained nonconvex
  quadratic programs with semidefinite-based finite branch-and-bound.
  \emph{Comput. Optim. Appl.} 43(2):181--195.

\bibitem[{Carlin et~al.(2007)Carlin, Lobo, \protect\BIBand{}
  Viswanathan}]{Carlin07}
Carlin BI, Lobo MS, Viswanathan S (2007) Episodic liquidity crises: Cooperative
  and predatory trading. \emph{J. Finance} 62(5):2235--2274.

\bibitem[{Chen \protect\BIBand{} Burer(2012)}]{Burer2012}
Chen J, Burer S (2012) Globally solving nonconvex quadratic programming
  problems via completely positive programming. \emph{Math. Program. Comput.}
  4:33--52.

\bibitem[{Chen et~al.(2015)Chen, Feng, \protect\BIBand{} Peng}]{Chen15}
Chen JN, Feng LM, Peng JM (2015) Optimal deleveraging with nonlinear temporary
  price impact. \emph{Eur. J. Oper. Res.} 244:240--247.

\bibitem[{Chen et~al.(2014)Chen, Feng, Peng, \protect\BIBand{} Ye}]{Chen14}
Chen JN, Feng LM, Peng JM, Ye YY (2014) Analytical results and efficient
  algorithm for optimal portfolio deleveraging with market impact. \emph{Oper.
  Res.} 62(1):195--206.

\bibitem[{Cont \protect\BIBand{} Wagalath(2016)}]{cont2016fire}
Cont R, Wagalath L (2016) Fire sales forensics: measuring endogenous risk.
  \emph{Math. Finance} 26(4):835--866.

\bibitem[{D'Hulster(2009)}]{d2009leverage}
D'Hulster K (2009) The leverage ratio: A new binding limit on banks .

\bibitem[{Fleming et~al.(1998)Fleming, Kirby, \protect\BIBand{}
  Ostdiek}]{fleming1998information}
Fleming J, Kirby C, Ostdiek B (1998) Information and volatility linkages in the
  stock, bond, and money markets. \emph{J. Financial Econom.} 49(1):111--137.

\bibitem[{{Gurobi Optimizer}(2020)}]{gurobi}
{Gurobi Optimizer} (2020) {G}urobi {I}nteractive {S}hell (win64), {V}ersion
  9.0.2 {C}opyright (c) 2020, {G}urobi {O}ptimization, {LLC} .

\bibitem[{Holthausen et~al.(1990)Holthausen, Leftwich, \protect\BIBand{}
  Mayers}]{Holthausen90}
Holthausen RW, Leftwich RW, Mayers D (1990) Large-block transactions, the speed
  of response, and temporary and permanent stock-price effects. \emph{J. Financial Econom.} 26(1):71--95.

\bibitem[{Hong et~al.(2011)Hong, Yang, \protect\BIBand{}
  Zhang}]{Hong11}
Hong  L, Yang Y, Zhang L (2011) Sequential convex approximations to joint chance constrained programs:
a Monte Carlo approach. \emph{Oper. Res.} 59: 617--630.


\bibitem[{Jiang et~al.(2018)Jiang, Li, \protect\BIBand{} Wu}]{Jiang2018}
Jiang RJ, Li D, Wu BY (2018) SOCP reformulation for the generalized trust
  region subproblem via a canonical form of two symmetric matrices. \emph{Math.
  Program. A} 169:531--563.

\bibitem[{Kyle \protect\BIBand{} Xiong(2001)}]{kyle2001contagion}
Kyle AS, Xiong W (2001) Contagion as a wealth effect. \emph{J. Finance}
  56(4):1401--1440.

\bibitem[{Li et~al.(2020)Li, Guo, Lai, Shi et~al.}]{lioptimal}
Li Y, Guo J, Lai KK, Shi J, et~al. (2020) Optimal portfolio liquidation with
  cross-price impacts on trading. \emph{Operational Research} 1--20.

\bibitem[{Linderoth(2005)}]{Linderoth05}
Linderoth J (2005) A simplicial branch-and-bound algorithm for solving
  quadratically constrained quadratic programs. \emph{Math. Program. B}
  103:251--282.

\bibitem[{Loridan(1982)}]{Lor82}
Loridan P (1982) Necessary conditions for $\epsilon$-optimality. \emph{Math.
  Program. Stud.} 19:140--152.

\bibitem[{Luo et~al.(2019)Luo, Bai, Lim, \protect\BIBand{} Peng}]{lblp19}
Luo HZ, Bai XD, Lim G, Peng JM (2019) New global algorithms for quadratic
  programming with a few negative eigenvalues based on alternative direction
  method and convex relaxation. \emph{Math. Program. Comput.} 11(1):119--171.

\bibitem[\protect\citeauthoryear{Luo, Ding, Peng, Jiang, and Li}{Luo
  et~al.}{2020}]{ldpjl20}
Luo  HZ, Ding XD, Peng JM, Jiang RJ,  Li D (2020) Complexity results and effective algorithms for the worst-case linear
  optimization under uncertainties.
\emph{INFORMS J. Comput.  https://doi.org/10.1287/ijoc.2019.0941\/}.


\bibitem[{Madhavan(2000)}]{Madhavan00}
Madhavan A (2000) Market microstructure: A survey. \emph{J. Financial Markets}
  3(3):205--258.


\bibitem[{McCormick(1976)}]{McCormick76}
McCormick GP (1976) Computability of global solutions to factorable nonconvex programs: Part II Convex
underestimating problems.  \emph{Math. Program.} 10:147--175.


\bibitem[{Mor\'{e}(1993)}]{More93}
Mor\'{e} JJ (1993) Generalization of the trust region problem. \emph{Optim.
  Methods Softw.} 2(1):189--209.

\bibitem[{Nesterov(1998)}]{Nesterov98}
Nesterov Y (1998) Semidefinite relaxation and nonconvex quadratic optimization.
  \emph{Optim. Methods Softw.} 9(1-3):141--160.

\bibitem[{Newcomb(1961)}]{Newcomb61}
Newcomb RW (1961) On the simultaneous diagonalization of two semi-definite
  matrices. \emph{Q. Appl. Math.} 19(2):144--146.

\bibitem[{Pasquariello \protect\BIBand{}
  Vega(2013)}]{pasquariello2013strategic}
Pasquariello P, Vega C (2013) Strategic cross-trading in the us stock market.
  \emph{Review of Finance} 19(1):229--282.


\bibitem[{Pham Dinh \protect\BIBand{}
  Le Thi(2014)}]{LeThi14}
Pham Dinh T and  Le Thi HA (2014)  Recent advances in {DC} programming and {DCA}, In Transactions on Computational Intelligence
XIII, 1--37. Springer, Berlin, Heidelberg.



\bibitem[{Polik \protect\BIBand{} Terlaky(2007)}]{Polik07}
Polik I, Terlaky T (2007) A survey of the {\it {S}}-lemma. \emph{SIAM Rev.}
  49(3):371--418.

\bibitem[{Sahinidis(1996)}]{int:Sahinidis96}
Sahinidis NV (1996) {BARON}: {A} general purpose global optimization software
  package. \emph{J. Global Optim.} 8:201--205.

\bibitem[{Sch\"{o}neborn \protect\BIBand{} Schied(2009)}]{Schied09b}
Sch\"{o}neborn T, Schied A (2009) Liquidation in the face of adversity: Stealth
  versus sunshine trading.
  \emph{http://papers.ssrn.com/sol3/papers.cfm?abstract.id=1007014} .

\bibitem[{Sias et~al.(2001)Sias, Starks, \protect\BIBand{} Titman}]{Sias01}
Sias RW, Starks LT, Titman S (2001) The price impact of institutional trading.
  \emph{http://papers.ssrn.com/sol3/papers.cfm?abstract.id=283779} .


\bibitem[{Saxena et~al.(2011)Saxena,  Bonami, \protect\BIBand{}
  Lee}]{Saxena11}
Saxena A,  Bonami P,   Lee J (2011) Convex relaxation of nonconvex mixed integer quadratically
  constrained programs: {P}rojected formulations.
\emph{ Math. Program.}  130:359--413.


\bibitem[{Tsoukalas et~al.(2019)Tsoukalas, Wang, \protect\BIBand{}
  Giesecke}]{Tsoukalas17}
Tsoukalas G, Wang J, Giesecke K (2019) Dynamic portfolio execution.
  \emph{Management Sci.} 65(5):2015--2040.

\bibitem[{Ye(1999)}]{Ye99}
Ye Y (1999) Approximating quadratic programming with bound and quadratic
  constraints. \emph{Math. Program.} 84(2):219--226.

\end{thebibliography}

\ECSwitch
\ECHead{Proofs of Statements}

In this e-companion, we provide the proofs of all the theorems, lemmas and propositions in the main body
of the paper.

\section{Proofs in Section \ref{sec:ODP}}
%\subsection{Proof of Proposition~\ref{prop:property-problem}.}
{\bf Proof of Proposition~\ref{prop:property-problem}.}  By Assumption~\ref{asmp1}, $l_1(-x_0)<0$ and hence $x_0^T\left(\Lambda+\frac{1}{2}\Gamma\right)x_0-e_0<0$.
Via a simple calculation, we  derive the following,
\begin{eqnarray*}
&& e_1(-x_0)=-x_0^T\left(\Lambda+\frac{1}{2}\Gamma\right)x_0+e_0>0,\\
&& \rho_1e_1(-x_0)-l_1(-x_0)=-(1+\rho_1)\left[x_0^T\left(\Lambda+\frac{1}{2}\Gamma\right)x_0-e_0\right]>0.
\end{eqnarray*}
The proof is completed. \Halmos\endproof

%\subsection{Proof of Proposition~\ref{prop-optimality}.}
{\bf Proof of Proposition~\ref{prop-optimality}.}
Since the feasible set of problem (\ref{OPDP}) is compact, there exists an optimal solution $y^*$. Let us denote $g_0(y)=\rho_1e_1(y)-l_1(y)$.
  Assume to the contrary that the leverage constraint is not active at the
optimal solution. Then $g_0(y^*)\neq0$. It is easy to see that the so-called
linear independence constraint qualification holds. Denote $y^*=(y^*_1,\ldots,y^*_m)^T$.
According to the first
order optimality condition, there exists $\mu^*=(\mu^*_0,\mu^*_1,\ldots,\mu^*_m,\mu^*_{m+1},\ldots,\mu^*_{2m})^T\ge0$ satisfying the following conditions
\begin{eqnarray}\label{kkt-0}
&&-\nabla e_1(y^*)+\mu^*_0\nabla g_0(y^*)+\sum_{i=1}^m\left(\mu^*_i\nabla g_i(y^*)+\mu^*_{m+i}\nabla g_{m+i}(y^*)\right)=0,\\\label{kkt-1}
&&\mu^*_0g_0(y^*)=0,\\\label{kkt-2}
&&\mu^*_ig_i(y^*)=0,\quad  i=1,\ldots,m,\\\label{kkt-3}
&&\mu^*_{m+i}g_{m+i}(y^*)=0,\quad i=1,\ldots,m,
 \end{eqnarray}
where
\[g_i(y)=y_i,\quad g_{m+i}(y)=-y_i-x_{0,i},\quad i=1,\ldots,m.\]
Since $g_0(y^*)\neq0$, we have $\mu^*_0=0$. Thus, equality (\ref{kkt-0}) becomes
\begin{eqnarray}\label{prop-eq-1}
  \sum_{j=1}^m\left[(\hat\lambda_{ij}-0.5\hat\gamma_{ij})y^*_j-\gamma_{ji}x_{0,j}\right]+\mu^*_i-\mu^*_{m+i}=0,~~i=1,\ldots,m,
 \end{eqnarray}
 where $\hat\lambda_{ij}=\lambda_{ij}+\lambda_{ji}$ and $\hat\gamma_{ij}=\gamma_{ij}+\gamma_{ji}$. Since matrices $\Gamma$  and $\Lambda$ are non-negative, we have  $\hat\lambda_{ij}\ge0$ and $\hat\gamma_{ij}\ge0$ for all $i,j$.
Since $x_0>0$, it is easy to see from (\ref{kkt-2}) and (\ref{kkt-3}) that $\mu^*_i$
and $\mu^*_{m+i}$, $1\le i \le m$, cannot be positive simultaneously. %By the assumption, both $\hat\Gamma$  and $\hat\Lambda$ are non-negative.
We consider the following cases:

Case (i): $\mu^*_i=0$, $\mu^*_{m+i}\neq 0$. From (\ref{kkt-3}), $y^*_i=-x_{0,i}$.
%By the assumption, $\gamma_{ij}-\gamma_{ji}\le 0$, $\forall i,j$.
Then, from $\mu^*\ge0$ and (\ref{prop-eq-1}), we obtain
\begin{eqnarray}\label{prop-eq-4}
0<\mu^*_{m+i}=\sum_{j=1}^m\left[(\hat\lambda_{ij}-0.5\hat\gamma_{ij})y^*_j-\gamma_{ji}x_{0,j}\right]\le\sum_{j=1}^m \hat\lambda_{ij}y^*_j+\frac{1}{2}\sum_{j=1}^m (\gamma_{ij}-\gamma_{ji})x_{0,j}\le 0,
 \end{eqnarray}
where the second inequality is due to $y^*\ge -x_0$ and $\hat\gamma_{ij}\ge0$ for all $i,j$, while the last inequality is due to $y^*\le 0$,   $\hat\lambda_{ij}\ge0$ and the assumption that $\Gamma$ is symmetric.   This is a contradiction.

Case (ii): $\mu^*_i=\mu^*_{m+i}=0$. Since $x_0>0$ and  $\Gamma$ is non-negative,
 it follows from  (\ref{prop-eq-1}) that
 \begin{eqnarray}\label{prop-eq-5}
   \sum_{j=1}^m (\hat\lambda_{ij}- 0.5\hat\gamma_{ij})y^*_j =\sum_{j=1}^mx_{0,j}\gamma_{ji}>0.
 \end{eqnarray}
 %from which  we must have $2\lambda_{ij_0}-\gamma_{ij_0}<0$ for some $j_0$ since $y^*\le 0$.
 Let us denote $J=\{j:\hat\lambda_{ij}-0.5\hat\gamma_{ij}<0,j=1,\ldots,m\}$.
 %If $J=\emptyset$, then $2\hat\lambda_{ij}-\hat\gamma_{ij}\ge0$ for $ j=1,\ldots,m$, which, by $y^*\le 0$ and (\ref{prop-eq-5}), implies that $0<\sum_{j=1}^m (2\hat\lambda_{ij}- \hat\gamma_{ij})y^*_j\le0$. This gives a contradiction. Thus, $J\neq\emptyset$.
 Since $y^*\le 0$ and by (\ref{prop-eq-5}), we must have $J\neq\emptyset$.
 Let us define
  \begin{eqnarray}\label{defn-j0}
 j_0\in\arg\max\limits_{j\in J}\{\hat\lambda_{ij}-0.5\hat\gamma_{ij}\}.
 \end{eqnarray}
 Note that $-x_0\le y^*\le 0$ implies that  $(\hat\lambda_{ij}-0.5\hat\gamma_{ij})y^*_j\le0$ for $j\not\in J$ and $(\hat\lambda_{ij}- 0.5\hat\gamma_{ij})y^*_j\le -(\hat\lambda_{ij}- 0.5\hat\gamma_{ij})x_{0,j}$ for $j\in J$.  It then follows from (\ref{prop-eq-5}) that
 \begin{eqnarray*}\label{prop-eq-6}
    \sum_{j=1}^mx_{0,j}\gamma_{ji}%=\sum_{j=1}^m (\hat\lambda_{ij}-0.5 \hat\gamma_{ij})y^*_j
   \le \sum_{j\in J} (\hat\lambda_{ij}- 0.5\hat\gamma_{ij})y^*_j\le (\hat\lambda_{ij_0}- 0.5\hat\gamma_{ij_0})y^*_{j_0} -\sum_{j\in J\backslash\{j_0\}}  (\hat\lambda_{ij}- 0.5\hat\gamma_{ij})x_{0,j},
 \end{eqnarray*}
 which in turn implies that
 \begin{eqnarray*}
  (\hat\lambda_{ij_0}- 0.5\hat\gamma_{ij_0})y^*_{j_0}&\ge&  \sum_{j=1}^m\gamma_{ji}x_{0,j}+\sum_{j\in J\backslash\{j_0\}}(\hat\lambda_{ij}- 0.5\hat\gamma_{ij})x_{0,j}\\
  &=& \sum_{j\not\in J}\gamma_{ji}x_{0,j}+\sum_{j\in J\backslash\{j_0\}}\left[\hat\lambda_{ij}x_{0,j}+\frac{1}{2}(\gamma_{ji}-\gamma_{ij})x_{0,j}\right]+x_{0,j_0}\gamma_{j_0i}\\
  &=&\sum_{j\not\in J}\gamma_{ji}x_{0,j}+\sum_{j\in J\backslash\{j_0\}}\hat\lambda_{ij}x_{0,j}+x_{0,j_0}\gamma_{j_0i},
 \end{eqnarray*}
 due to the assumption that $\Gamma$ is symmetric.
Since $\hat\lambda_{ij_0}-0.5 \hat\gamma_{ij_0}<0$ by $j_0\in J$, the above inequality yields
 \begin{eqnarray}\label{prop-eq-7}
 y^*_{j_0}\le \frac{\sum_{j\not\in J}\gamma_{ji}x_{0,j} +\sum_{j\in J\backslash\{j_0\}}\hat\lambda_{ij}x_{0,j}+x_{0,j_0}\gamma_{j_0i}}{\hat\lambda_{ij_0}- 0.5\hat\gamma_{ij_0}} .
\end{eqnarray}
Combining (\ref{prop-eq-7}) with $y^*_{j_0}\ge -x_{0,j_0}$ and $\gamma_{j_0i}=\gamma_{ij_0}$, we derive
\[0\le y^*_{j_0}+x_{0,{j_0}} \le  \frac{\sum_{j\not\in J}x_{0,j}\gamma_{ji} +\sum_{j\in J}\hat\lambda_{ij}x_{0,j}}{\hat\lambda_{ij_0}-0.5 \hat\gamma_{ij_0}}
 <0,\]
 which gives a contradiction.

Since for any $1\le i\le m$, the above cases are not possible,
we must have $\mu^*_i\neq0$, $\mu^*_{m+i}= 0$ for all $1\le i\le m$. Then from (\ref{kkt-2}), $y^*_i=0$, $i=1,\ldots,m$. Thus
$g_0(y^*)\le 0$ becomes $ l_0-\rho_1 e_0\le 0$,
which contradicts the assumption that $ l_0-\rho_1 e_0> 0$.  This completes the proof of the proposition.\Halmos\endproof

%\subsection{Proof of Proposition~\ref{prop-priority}. }
{\bf Proof of Proposition~\ref{prop-priority}. }
Assume to the contrary that $y^*$ is the optimal solution and $-x_{0,j}\le y_{j}^*< y_{i}^*\le0$. We will find a direction
$d\in{\mathbb R}^m$ along which the liability $l_1$ is strictly decreasing but the equity
$e_1$ is strictly increasing. Let $d=(d_1,\ldots,d_m)^T$, where $d_i=-1$, $d_j=1$, and $d_k=0$, $\forall k\neq i,j$.  By
the expressions of $l_1(y)$ and $e_1(y)$, we have
\begin{eqnarray}  \label{prop-priority-eq1}
 \nabla l_1(y)=(2\hat\Lambda+\hat\Gamma)y+p_0,\quad\quad
 \nabla e_1(y)=-(2\hat\Lambda-\hat\Gamma)y+\Gamma^Tx_0.
\end{eqnarray}
Because $p_{0,i}=p_{0,j}$ and $x_{0,i}=x_{0,j}$, by conditions (i)-(ii) and the fact that $y_j^*-y_i^*<0$, we can derive from \eqref{prop-priority-eq1} that
\begin{eqnarray*} \nonumber
 \nabla l_1(y^*)^Td&=&\sum_{k=1}^m(2\hat\lambda_{jk}+\hat\gamma_{jk}-2\hat\lambda_{ik}-\hat\gamma_{ik})y_k^*\\
 &\le&[2(\hat\lambda_{jj}-\hat\lambda_{ij})+\hat\gamma_{jj}-\hat\gamma_{ij}]y_j^*+[2(\hat\lambda_{ji}-\hat\lambda_{ii})+\hat\gamma_{ji}-\hat\gamma_{ii}]y_i^*\\
 &\le& [2(\hat\lambda_{jj}-\hat\lambda_{ij})+\hat\gamma_{jj}-\hat\gamma_{ij}](y_j^*-y_i^*)<0.  \label{prop2-Liquidation-ineq1}
\end{eqnarray*}
Meanwhile, the equation \eqref{prop-priority-eq1} implies that
\begin{eqnarray}
 \nabla e_1^T(y^*)d=\sum_{k=1}^m\left[(-2\hat\lambda_{jk}+\hat\gamma_{jk}+2\hat\lambda_{ik}-\hat\gamma_{ik})y^*_k+(\gamma_{kj}-\gamma_{ki})x_{0,k}\right]. \label{prop2-Liquidation-ineq3}
\end{eqnarray}
Because $-x_{0,k}\le y^*_k\le0$ holds for all $k$,  we follow from condition (i), (ii) and (iv) that % for all $k\neq i,j$,
\begin{eqnarray} \nonumber
&& (-2\hat\lambda_{jk}+\hat\gamma_{jk}+2\hat\lambda_{ik}-\hat\gamma_{ik})y^*_k+(\gamma_{kj}-\gamma_{ki})x_{0,k}\\\nonumber
&&=2(\hat\lambda_{ik}-\hat\lambda_{jk})y^*_k+(\hat\gamma_{jk}-\hat\gamma_{ik})y^*_k+(\gamma_{kj}-\gamma_{ki})x_{0,k}\\\nonumber
&&\ge -(\hat\gamma_{jk}-\hat\gamma_{ik})x_{0,k}+(\gamma_{kj}-\gamma_{ki})x_{0,k}\\
&&=\frac{1}{2}[(\gamma_{kj}-\gamma_{jk})+(\gamma_{ik}-\gamma_{ki})]x_{0,k}\ge0,\quad \forall k\neq i,j.\label{prop2-Liquidation-ineq4}
 \end{eqnarray}
Thus, according to the inequality (\ref{prop2-Liquidation-ineq4}), we can relax the gradient (\ref{prop2-Liquidation-ineq3}) to
\begin{eqnarray*} \nonumber
 \nabla e_1^T(y^*)d&\ge& [2(\hat\lambda_{ij}-\hat\lambda_{jj})+\hat\gamma_{jj}-\hat\gamma_{ij}]y^*_j+(\gamma_{jj}-\gamma_{ji})x_{0,j}\nonumber\\
 &&+[2(\hat\lambda_{ii}-\hat\lambda_{ji})+\hat\gamma_{ji}-\hat\gamma_{ii}]y^*_i+(\gamma_{ij}-\gamma_{ii})x_{0,i} \nonumber\\
 &=& 2(\hat\lambda_{ij}-\hat\lambda_{jj})y^*_j+(\hat\gamma_{jj}-\hat\gamma_{ij})y^*_j
   +2(\hat\lambda_{ii}-\hat\lambda_{ji})y^*_i\nonumber\\
 &&  +(\hat\gamma_{ji}-\hat\gamma_{ii})y^*_i+(\gamma_{jj}+\gamma_{ij}-\gamma_{ji}-\gamma_{ii})x_{0,i} \nonumber\\
&\ge& 2(\hat\lambda_{ij}-\hat\lambda_{jj})(y_j^*-y^*_i)
+[-(\hat\gamma_{jj}-\hat\gamma_{ij})-(\hat\gamma_{ji}-\hat\gamma_{ii})+(\gamma_{jj}-\gamma_{ii})]x_{0,i} \\
%&=& 2(\hat\lambda_{ij}-\hat\lambda_{jj})(y_j^*-y^*_i)+[-\hat\gamma_{jj}+\hat\gamma_{ii}+\gamma_{jj}-\gamma_{ii}]x_{0,i} \\
&>&0, \label{prop2-Liquidation-ineq2}
\end{eqnarray*}
where the second inequality follows from conditions (i)-(iii) and the fact $-x_{0,i}=-x_{0,j}\le y_{j}^*< y_{i}^*\le0$, while the last inequality is due to condition (i) and the fact that $y_{j}^*< y_{i}^*$.

 Therefore, there exist a sufficiently small $\delta>0$ such that $-x_0\le y^*+\delta d \le0$ and the trading strategy
$y^*+\delta d$ leads to strictly larger equity and smaller liability after trading. This contradicts the optimality of $y^*$.
\Halmos\endproof

%\subsection{ The proof of Lemma~\ref{lem-SD}.}
%We first recall a well-known result  regarding the simultaneous diagonalizability of two semi-definite matrices.
%The $n\times n$ real symmetric matrices $A$ and $B$ are called simultaneously diagonalizable (SD) if there exists a nonsingular matrix $D$ such that both $D^T AD$ and $D^T BD$ are diagonal.
%For $d\in{\mathbb R}^n$, let ${\rm diag}(d)$ denote  an $n\times n$ diagonal matrix with its diagonal equal to $d$. The following result is from  [\cite{Newcomb61}, Theorem].
%
%
%\begin{lemma}\label{lem-SD}
%Let $A$ and $B$ be the $n\times n$ real symmetric positive semi-definite matrices, and let $s={\rm rank} (B)$, $q={\rm rank} (A)$  and $r={\rm rank}(A+B)$. Then there exist  a nonsingular matrix $D$   such that
%\begin{eqnarray*}%\label{SD}
%D^TBD={\rm diag}(\delta_1,\ldots,\delta_s,0,\ldots,0),~~D^TAD={\rm diag}(\theta_1,\ldots,\theta_r,0,\ldots,0),
%\end{eqnarray*}
% where  $s<r$,  $0<\delta_i\le 1$ for $i=1,\ldots,s$,   $\theta_i=1-\delta_i$ for $i=1,\ldots,s$ and $\theta_{i}=1$ for $i=s+1,\ldots,r$, and there are $q$ nonzero numbers in $\theta_1,\ldots,\theta_r$.
%\end{lemma}

{\bf The proof of Lemma~\ref{lem-SD}.}
In order to find a nonsingular matrix $D$ in the simultaneous diagonalizability  method, we give the proof of Lemma~\ref{lem-SD} slightly different from that in \cite{Newcomb61} as follows.

%{\bf The proof of Lemma~\ref{lem-SD}}~~
Let $A$ and $B$ be two $n\times n$ real symmetric  positive semidefinite matrices, and let  $r={\rm rank} (A+B)$ and $s={\rm rank} (B)$. Then there exists a nonsingular matrix $Q\in{\mathbb R}^{n\times n}$ such that $Q^T(A+B)Q=\left(\begin{array}{cc}I_r&0\\0&0\end{array}\right)$, where $I_r$ denotes the unit matrix of order $r$. Let $Q^TBQ=\left(\begin{array}{cc}G_{11}&G_{12}\\G_{12}^T&G_{22}\end{array}\right)$ with $G_{11}\in{\mathbb R}^{r\times r}$. Since $A$ is positive semidefinite, we have $Q^T(A+B)Q\succeq Q^TBQ$ and hence $\left(\begin{array}{cc}I_r&0\\0&0\end{array}\right)\succeq \left(\begin{array}{cc}G_{11}&G_{12}\\G_{12}^T&G_{22}\end{array}\right)$, from which we can infer that
%$I_r\succeq B_{11}\succeq 0$, $G_{12}=0$ and $G_{22}=0$. Thus,
 $Q^TBQ=\left(\begin{array}{cc}G_{11}&0\\0&0\end{array}\right)$ and $I_r\succeq G_{11}\succeq 0$. Note that ${\rm rank} (G_{11})={\rm rank} (B)=s$. Let $\delta_1,\ldots, \delta_s$ be the positive eigenvalues of matrix $G_{11}$. Since $I_r\succeq G_{11}\succeq 0$, we obtain $s< r$ and $0< \delta_i\le 1$ for $i=1,\ldots,s$.    Since $G_{11}$ is a real symmetric matrix of order $r$, there exists an orthogonal matrix $D_r\in{\mathbb R}^{r\times r}$ such that $D_r^TG_{11}D_r={\rm diag}(\delta_1,\ldots,\delta_s,0,\ldots,0)$.   Let us choose $D=Q\left(\begin{array}{cc}D_r&0\\0&I_{n-r}\end{array}\right)$. It is then easy to check that $D$ is nonsingular, $D^TBD={\rm diag}(\delta_1,\ldots,\delta_s,0,\ldots,0)$ and
$D^TAD={\rm diag}(\theta_1,\ldots,\theta_r,0,\ldots,0)$, where $\theta_i=1-\delta_i$ for $i=1,\ldots,s$ and $\theta_{i}=1$ for $i=s+1,\ldots,r$. Furthermore,  since $D$ is nonsingular and $q={\rm rank} (A)$, there are $q$ nonzero numbers in $\theta_1,\ldots,\theta_r$. \Halmos\endproof

%{\bf The proof of Proposition~\ref{prop-SD}}. Since $B^-$ and $A^-$ are symmetric positive semi-definite matrices with ${\rm rank} (B^-)=s$ and ${\rm rank} (A^-)=q$, by Lemma~\ref{lem-SD}, $B^-$ and $A^-$ are simultaneously diagonalizable with the nonsingular matrix $D\in{\mathbb R}^{m\times m}$, i.e.,
%\begin{eqnarray}
%%\left\{\begin{array}{l}
% D^T B^-D ={\rm diag}(\delta_1,\ldots,\delta_s,0,\ldots,0),\quad\quad D^T A^-D ={\rm diag}(\theta_1,\ldots,\theta_r,0,\ldots,0),
%%\end{array}\right.
%\end{eqnarray}
%where $s<r$, $0<\delta_i\le 1$ for $i=1,\ldots,s$,  $\theta_i=1-\delta_i$ for $i=1,\ldots,s$  and $\theta_{i}=1$ for  $i=s+1,\ldots,r$, and there are $q$ nonzero numbers in $\theta_1,\ldots,\theta_r$. With the facts in (\ref{B-decomp})  and (\ref{SD}) and the transformation  $y=Dz$, we can transform the original problem \eqref{OPDP-sym} to \eqref{OPDP-SD}.
%
%Moreover, these two problems~(\ref{OPDP-sym}) and (\ref{OPDP-SD})  are equivalent in the sense that they have the same optimal value and, if $z^*$ is a  global optimal solution of Problem (\ref{OPDP-SD}), then $y^*=Dz^*$ is a global optimal solution of Problem~(\ref{OPDP-sym}). \Halmos\endproof

\section{Proofs in Section~\ref{sec:SCO}}

%\subsection{Proof of Lemma~\ref{lem-F-nonemptiness}.}

{\bf Proof of Lemma~\ref{lem-F-nonemptiness}.}
i) We first observe that
 \begin{eqnarray}\label{eq-g}
    \hat{g}_{\xi}(z)=\hat{g}(z)+\sum_{i=1}^{r}\theta_i(z_i-\xi_i)^2+\rho_1\sum_{i=1}^{s}\delta_i(z_i-\xi_i)^2.
 \end{eqnarray}
Since $\bar z\in\hat{\cal F}$ and $\bar\xi_i=\bar z_i$ for $i=1,\ldots,r$, we can infer from (\ref{eq-g}) that $\bar y\in \hat{\cal F}_{\bar\xi}$  and hence $\hat{\cal F}_{\bar\xi}\neq\emptyset$.
It is also easy to see that $\hat{\cal F}_{\bar\xi}$ is a   closed convex set. In addition,
 we  can follow from (\ref{eq-g}) that $\hat{g}_{\xi}(z)\ge \hat{g}(z)$ and hence $\hat{\cal F}_{\bar\xi}\subseteq\hat{\cal F}$.

ii) By contradiction, suppose that  $ {\rm int}\hat{\cal F}_{\bar\xi}=\emptyset$. This means that
\[\hat{g}_{\bar\xi}(z)\ge 0,\quad \forall  z\in {\cal Z}.\]
Since $\bar z\in\hat{\cal F}$ and $\bar\xi_i=\bar z_i$ for $i=1,\ldots,r$, we can deduce from (\ref{eq-g}) that
$0\le \hat{g}_{\bar\xi}(\bar z)=\hat{g}(\bar z) \le0$,
which in turn  implies that  $\hat{g}(\bar z)=0$ and
\begin{eqnarray*}\label{lem-slater-prob}
 \bar z\in \arg\min_{z \in {\cal Z}}\hat{g}_{\bar\xi}(z) .
\end{eqnarray*}
From the necessary optimality condition, there exist  $\eta$ and $\mu\in{\mathbb R}^m_+$ such that
\begin{eqnarray}\label{lem-prob-kkt}
  \nabla \hat{g}_{\bar\xi}(\bar z)+D(\eta-\mu)  =0,\quad \eta^TD\bar{z}=0,  \quad \mu^T(D\bar{z}+x_0)=0 .
\end{eqnarray}
Since $\bar\xi_i=\bar z_i$ for $i=1,\ldots,r$, from (\ref{eq-g}), we have $\nabla \hat{g}_{\bar\xi}(\bar z)=\nabla \hat{g}(\bar z)$. Let $\bar y=D\bar z$. Clearly, $-x_0\le\bar y\le 0$ and $g(\bar y)=\hat{g}(\bar z)=0$, so $\bar y$ is a feasible solution to  problem (\ref{OPDP-sym}). We also follow that  $\nabla \hat{g}(\bar z)=P\nabla g(\bar z)$. Note that $P$ is nonsingular. It then follows from (\ref{lem-prob-kkt}) that
\begin{eqnarray}\label{lem-prob-kkt2}
  \nabla  g(\bar z)+ (\eta-\mu)  =0,\quad \eta^T \bar{y}=0,  \quad \mu^T(\bar{y}+x_0)=0 .
\end{eqnarray}
Let us define $I_1=\{i\mid \bar y_i=0,~i=1,\ldots,m\}$, and $I_2=\{i\mid \bar y_i=-x_{0,i},~i=1,\ldots,m\}$.
It then follows from   (\ref{lem-prob-kkt2}) that
\begin{eqnarray*}\label{lem-prob-kkt3}
&& \nabla g(\bar y)+\sum_{i\in I_1} \eta_ie_i-\sum_{i\in I_2} \mu_ie_i =0,
\end{eqnarray*}
where $e_i$ is the $i$th unit vector of ${\mathbb R}^m$. This, together with the feasibility of $\bar y$, contradicts the fact that the gradients $\nabla g(\bar y)$, $e_i$, $i\in I_1\cup I_2$ are linear independent. \Halmos\endproof

%\subsection{Proof of Lemma~\ref{lem1-SCA}.}
{\bf Proof of Lemma~\ref{lem1-SCA}.}
i)  Since $z^0\in\hat{\cal F}$ and $\xi^0=(\xi_1^0,\ldots,\xi_r^0)^T$ with $\xi_i^0=z_i^0$ for $i=1,\ldots,r$, we obtain from (\ref{eq-g}) that $z^0\in\hat{\cal F}_{\xi^0}$ and hence $\hat{\cal F}_{\xi^0}\neq\emptyset$. We see from Step 1 that $z^1\in\hat{\cal F}_{\xi^0}$. It then follows from Lemma~\ref{lem-F-nonemptiness} that  $z^1\in\hat{\cal F}$.  By induction,   we  conclude that  $\{z^k\}\subseteq\hat{\cal F}$.

ii) From claim (i), we have  $\{z^k \}\subseteq\hat{\cal F}$. Note from Step 1 that $\xi^{k}=(\xi_1^{k},\ldots,\xi_r^{k})^T$ with $\xi_i^{k}=z_i^{k}$ for $i=1,\ldots,r$. It is  easy to check from (\ref{eq-g}) that $z^k \in\hat{\cal F}_{\xi^k}$ for all $k$. Note that
 \begin{eqnarray}\label{eq-f}
    \hat{f}_{\xi}(z)=\hat{f}(z)+ \sum_{i=1}^{s}\delta_i(z_i-\xi_i)^2.
 \end{eqnarray}
From Step 1, since $z^{k+1}$ is the optimal solution of problem~(\ref{ODP-linearappro}) with  $\xi=\xi^{k}$ and
  by $z^k \in\hat{\cal F}_{\xi^k}$, we can deduce  that
\begin{eqnarray*}\label{prop-algo1-ineq1}
\hat{f}(z^{k+1})+\sum_{i=1}^{s}\delta_i(z_i^{k+1}-\xi^k_i)^2 \le \hat{f}(z^{k})+\sum_{i=1}^{s}\delta_i(z_i^{k}-\xi^k_i)^2,
\end{eqnarray*}
 which,  together with the fact that $\xi_i^{k}=z_i^{k}$ for $i=1,\ldots,r$ and $s<r$, yields
\[ \hat{f}(z^{k})-\hat{f}(z^{k+1}) \ge \sum_{i=1}^{s}\delta_i(\xi_i^{k+1}-\xi^k_i)^2, ~~~~~\forall k\ge 1.\]
The proof is finished.\Halmos\endproof

%\subsection{Proof of  Lemma~\ref{lem3-SCA}.}

{\bf Proof of  Lemma~\ref{lem3-SCA}.}
By Lemma~\ref{lem1-SCA}, we have that
\begin{eqnarray}\label{lem1-algo1-ineq1}
\hat{f}(z^{k})-\hat{f}(z^{k+1}) \ge \sum_{i=1}^{s}\delta_i(\xi^{k+1}_i-\xi^k_i)^2
\end{eqnarray}
for all $k$. Thus,
 $\{\hat{f}(z^k )\}$ is a non-increasing sequence. By Lemma~\ref{lem1-SCA}, $\{z^k \}\subseteq\hat{\cal F}$. Since $\hat{\cal F}$ is compact, $\{\hat{f}(z^k)\}$ is bounded. Hence $\{\hat{f}(z^k )\}$ converges and consequently
 \begin{eqnarray}\label{algo1-fk-lim}
 \hat{f}(z^{k})-\hat{f}(z^{k+1})\to 0~(k\to\infty).
 \end{eqnarray}
By the fact that $\delta_i>0$ for $i=1,\ldots,s$, it then follows  from (\ref{lem1-algo1-ineq1}) that
\begin{eqnarray}\label{algo1-tk-lim}
\lim_{k\to \infty} (\xi^{k+1}_i-\xi^k_i)=0,\quad i=1,\ldots,s.
\end{eqnarray}

Now let $(\bar z,\bar\xi)$ be an accumulation point of the sequence $\{(z^k,\xi^k)\}$, and let ${\cal K}\subset\{1,2,\ldots\}$ be such that $\{(z^k,\xi^k)\}_{\cal K}\to (\bar z,\bar\xi)$. The closedness of $\hat{\cal F}$ and $\{z^k \}\subseteq\hat{\cal F}$ imply $\bar z\in\hat{\cal F}$. Note from Step 1 that $\xi^k_i=z^k_i$, $i=1,\ldots,r$ for all $k$, so $\bar\xi_i= \bar z_i$ for $i=1,\ldots,r$. We then follow from (\ref{eq-g}) that $\bar z\in\hat{\cal F}_{\bar\xi}$.
In the following, we will prove that for any $z\in\hat{\cal F}_{\bar\xi}$, it holds that
\begin{eqnarray}\label{thm2-ineq-opt}
\hat{f}_{\bar\xi}(z)\ge \hat{f}_{\bar\xi}(\bar{z}).
\end{eqnarray}
%This, together with $\bar z\in\hat{\cal F}_{\bar\xi}$,  implies that $\bar z$ is the optimal solution of problem~(\ref{ODP-linearappro}) with  $\xi=\bar\xi$.
% As a result,  by Proposition~\ref{prop2-optimality}, $\bar z$  is  a KKT point of  Problem (\ref{OPDP-SD}).
%
Indeed, for any $z\in\hat{\cal F}_{\bar\xi}$, we consider  two cases as follows.

Case (i): $z\in\hat{\cal F}_{\xi^k}$ for a sufficiently large $k\in{\cal K}$. In this case, since $z^{k+1}$ is an optimal solution of  problem~(\ref{ODP-linearappro}) with  $\xi=\xi^{k}$, we have
 \begin{eqnarray}\label{thm2-ineq2}
 \hat{f}_{\xi^k}(z)\ge \hat{f}_{\xi^k}(z^{k+1}).
\end{eqnarray}
Since $\xi^k_i=z^k_i$, $i=1,\ldots,r$ for all $k$ and $s<r$, we follow from  (\ref{eq-f}) that
 \begin{eqnarray}\label{fk-eq}
\hat{f}_{\xi^k}(z^{k+1})=\hat{f}(z^{k+1})+ \sum_{i=1}^{s}\delta_i(\xi^{k+1}_i-\xi^k_i)^2.
\end{eqnarray}
Note from (\ref{algo1-fk-lim}) that $\lim\limits_{k\to\infty,k\in{\cal K}} \hat{f}(z^{k+1})=\hat{f}(\bar z)$. Note from (\ref{eq-f}) and $\bar\xi_i=\bar z_i$, $i=1,\ldots,r$ that $\hat{f}(\bar{z})= \hat{f}_{\bar\xi}(\bar{z})$. It then follows from (\ref{fk-eq}) and (\ref{algo1-tk-lim}) that
 \begin{eqnarray}\label{lim-fk}
\lim\limits_{k\to\infty,k\in{\cal K}}\hat{f}_{\xi^k}(z^{k+1})=\hat{f}_{\bar\xi}(\bar{z}).
\end{eqnarray}
Also, $\lim\limits_{k\to\infty,k\in{\cal K}}\hat{f}_{\xi^{k}}(z)=\hat{f}_{\bar\xi}(z)$. Then taking the limit in (\ref{thm2-ineq2}) as $k\to\infty$ and $k\in{\cal K}$  yields  (\ref{thm2-ineq-opt}).

Case (ii): $z\not\in\hat{\cal F}_{\xi^k}$ for a sufficiently large $k\in{\cal K}$. Since $\bar{z}\in\hat{\cal F}$ and $\bar\xi_i=\bar{z}_i$, $i=1,\ldots,r$,
by Lemma~\ref{lem-F-nonemptiness}, ${\rm int}{\cal F}_{\bar\xi}\neq\emptyset$.
 Then there exists $\hat{z}\in{\cal Z}$ such that $\hat{g}_{\bar\xi}(\hat{z})<0$. Let $\delta=-\hat{g}_{\bar\xi}(\hat{z})>0$. Since $\hat{g}_{\xi}(\hat{z})$ is continuous in $\xi$ and by $\{\xi^k\}_{\cal K}\to\bar\xi$, we have that for sufficiently large $k\in\cal K$, one has
\[\hat{g}_{\xi^k}(\hat{z})\le -\delta/2<0.\]
For given  $z\in\hat{\cal F}_{\bar\xi}\setminus{\cal F}_{\xi^k}$,
let us define
$\rho_k=\max\left\{0,\hat{g}_{\xi^k}(z)\right\}$.
Clearly,  $\rho_k>0$. Since $\{\xi^k\}_{\cal K}\to\bar\xi$ and $z\in\hat{\cal F}_{\bar\xi}$, it is easy to check that
\begin{eqnarray}\label{limit:rho}
\lim\limits_{k\to\infty,k\in{\cal K}}\rho_k=0.
\end{eqnarray}
Let us define
\begin{eqnarray}   \label{defn-lambda}
\lambda_k= 2\rho_k/(2\rho_k+\delta).
\end{eqnarray}
Clearly, $0<\lambda_k<1$. Let us choose $\hat{z}^k=(1-\lambda_k)z+\lambda_k\hat{z}$.
Since $\hat{g}_{\xi^k}(z)$ is a convex function in $z$ and by (\ref{defn-lambda}), we have
\begin{eqnarray*}   %\label{lem3-ineq1}
\hat{g}_{\xi^k}(\hat{z}^k) \le (1-\lambda_k)\hat{g}_{\xi^k}(z)+\lambda_k\hat{g}_{\xi^k}(\hat{z})
 \le  (1-\lambda_k)\rho_k+\lambda_k(-\delta/2)=0,
\end{eqnarray*}
which in turn implies that $\hat{z}^k\in\hat{\cal F}_{\xi^k}$.
Therefore, we have
\begin{eqnarray*} \label{lem3-ineq2}
{\rm dist}(z,\hat{\cal F}_{\xi^k}):=\min\{\|w-z\|:w\in\hat{\cal F}_{\xi^k}\}\le \|\hat{z}^k-z\|.
\end{eqnarray*}
%where ${\rm dist}(x,{\cal F}_{\xi^k})$ denotes the distance between $x\in [-x_0,0]$ and ${\cal F}_{t^k}$.
Note that
$\hat{z}^k-z=\lambda_k(\hat{z}-z)$. Since $\hat{\cal F}_{\xi^k}$ is a nonempty closed convex set, by the projection theorem, there exists a unique $\tilde{z}^k\in\hat{\cal F}_{\xi^k}$ such that
\begin{eqnarray*} %\label{lem3-ineq3}
\|\tilde{z}^k-z\|={\rm dist}(z,\hat{\cal F}_{\xi^k})\le \|\hat{z}^k-z\|=\lambda_k\|\hat{z}-z\|\le 2\rho_k\delta^{-1}\|\hat{z}-z\|.
\end{eqnarray*}
From the above relation and (\ref{limit:rho}), we obtain
\begin{eqnarray}\label{limit:yk}
\lim_{k\in\cal K\to \infty} \tilde{z}^k=z.
\end{eqnarray}
 Since $\tilde{z}^k\in\hat{\cal F}_{\xi^k}$ and $z^{k+1}$ is an optimal solution of  problem~(\ref{ODP-linearappro}) with  $\xi=\xi^{k}$, we have
 \begin{eqnarray}\label{thm2-ineq3}
 \hat{f}_{\xi^k}(z^{k+1})\le  \hat{f}_{\xi^k}(\tilde{z}^k).
\end{eqnarray}
By (\ref{lim-fk}),  $\lim\limits_{k\in{\cal K}\to\infty} \hat{f}_{\xi^k}(z^{k+1})= \hat{f}_{\bar\xi}(\bar{z})$.
By  (\ref{limit:yk}) and $\{\xi^k\}_{\cal K}\to \bar\xi$, we obtain  $\lim\limits_{k\in{\cal K}\to\infty} \hat{f}_{\xi^k}(\tilde{z}^k)= \hat{f}_{\bar\xi}(z)$. Taking the limit in (\ref{thm2-ineq3}) as $k\to\infty$ and $k\in{\cal K}$ then gives rise to (\ref{thm2-ineq-opt}).
The proof of the theorem is completed.
\Halmos\endproof

\section{Proofs in Section~\ref{sect:GOA}}
%\subsection{Proof of Theorem \ref{thm1-CR}.}
{\bf Proof of Theorem \ref{thm1-CR}.}
 %Since problem  (\ref{ODP:Bounded-CR}) is a convex relaxation of problem (\ref{OPD:Bounded}),
We first have $\hat{f}^{*}_{[l,u]}\ge \hat{v}^*_{[l,u]}(\bar{y})$. %Note that $0<\delta_i\le1$ for $i=1,\ldots,s$.
  It follows  that
\begin{eqnarray*} \nonumber
  \hat{f}(\bar z)-\hat{f}^*_{[l,u]}&\le& \hat{f}(\bar z)-\hat{v}^*_{[l,u]}=\sum_{i=1 }^{s}\delta_i(\bar t_i-  \bar z_i^2) \\
    & \le&\sum_{i=1 }^{s} (\bar t_i-  \bar z_i^2) \label{thm2-CR-ineq}\\\nonumber
  &\le& \sum_{i=1 }^{s} [-\bar z_i^2+(l_i+u_i)\bar z_i -l_iu_i] \\\nonumber
 &\le& \frac{1}{4}\sum_{i=1 }^{s}(u_i-l_i)^2\le \frac{s}{4}\|u-l\|^2_\infty,
  \end{eqnarray*}
 where the second and third  inequalities follow  from the fact that $\delta_i\in(0,1]$ for $i=1,\ldots,s$ and the constraints on $(t_i, z_i)$ in (\ref{ODP:Bounded-CR}).

 Similarly, using the feasibility of $(\bar{z},\bar{t})$ and the fact that $0<\theta_i\le1$ and $0<\delta_i\le1$ for each $i$, we can derive that
\begin{eqnarray*} \nonumber
\hat{g}(\bar z)&=&\psi(\bar z)-\sum_{i=1}^{r}\theta_i\bar t_i-\rho_1\sum_{i=1}^{s}\delta_i\bar t_i+\sum_{i=1}^{r}\theta_i(\bar t_i-\bar z_i^2) +\rho_1\sum_{i=1}^{s}\delta_i(\bar t_i-\bar z_i^2)\\
&\le& \sum_{i=1}^{r}(\bar t_i-\bar z_i^2) +\rho_1\sum_{i=1}^{s}(\bar t_i-\bar z_i^2)\\
&\le&(r+\rho_1s)\max_{i=1,\ldots,r}\{\bar t_i-  \bar z_i^2\}\le\frac{r+\rho_1s}{4}\|u-l\|^2_\infty.
  \end{eqnarray*}
This completes the proof.  \Halmos\endproof

%\subsection{Proof of Lemma~\ref{lemma-SCOBBA}.}

{\bf Proof of Lemma~\ref{lemma-SCOBBA}.}
  Recall that $(z^k,t^k)$ is the optimal solution of problem~(\ref{ODP:Bounded-CR}) over $\Delta^k$. From the proof of Theorem~\ref{thm1-CR}, we can obtain Statement (i).
 Now we prove Statement (ii).
 Since $(z^k,t^k)$ is the optimal solution of problem~(\ref{ODP:Bounded-CR}) over $\Delta^k$ and $v^k$ is its optimal value, by Theorem~\ref{thm1-CR}, we can infer that $\hat{f}(z^k)-v^k\le\epsilon$ and $\hat{g}(z^k)\le\epsilon$. Thus, $z^k$ is $\epsilon$-feasible for Problem (\ref{OPDP-SD}).
From Steps 2 and (S3.5)  of Algorithm~\ref{SCOBBA} , we see that $\hat{f}(z^k)\ge v^*=\hat{f}(z^*)$ for all $k$. Thus,
\begin{eqnarray}\label{lemma-SCOBBA-ineq1}
v^k- v^*\ge v^k-\hat{f}(z^k)\ge-\epsilon,
\end{eqnarray}
 so the stopping criterion of the algorithm is satisfied and then the algorithm   stops. Let $\hat{f}^*$ denotes the optimal value of Problem (\ref{OPDP-SD}). By Step (S3.1), $v^k$ is the smallest lower bound. Thus $\hat{f}^*\ge v^k$. It then  follows from (\ref{lemma-SCOBBA-ineq1}) that
 \begin{eqnarray*}\label{lemma-SCOBBA-ineq2}
    \hat{f}(z^*)=v^*\le v^k+\epsilon\le \hat{f}^*+\epsilon.
\end{eqnarray*}
This, together with   $\hat{g}(z^*)\le\epsilon$, implies that   $z^*$ is a global  $\epsilon$-solution to problem (\ref{OPDP-SD}). \Halmos\endproof

%\subsection{Proof of Theorem~\ref{thm1-algo1-complexity}.}
%The proof is similar to the proof of Theorems 3.3 in \cite{lblp19}, obtained by using Lemmas~\ref{lemma1-SCOBBA} and \ref{lemma-SCOBBA-stop}.
%\proof{Proof.}
{\bf Proof of Theorem~\ref{thm1-algo1-complexity}.}
At the $k$-th iteration, if the selected node $[\Delta^k, v^k, (z^k,t^k)]$ with the smallest lower bound $v^k$ satisfies either $u_{i^*}^k- l_{i^*}^k\le \frac{2\sqrt{\epsilon}}{\sqrt{r+\rho_1s}}$ for the chosen $i^*$ in partition  or
$\|u^k- l^k\|_\infty\le \frac{2\sqrt{\epsilon}}{\sqrt{r+\rho_1s}}$,
then by Lemma~\ref{lemma-SCOBBA}, the algorithm  stops and yields an  $\epsilon$-optimal solution $z^*$ to Problem (\ref{OPDP-SD}).
At the $k$-th iteration, if the algorithm  does not stop, then by Lemma~\ref{lemma-SCOBBA}(ii), $s_{i^*}^k- (t_{i^*}^k)^2> \frac{\epsilon}{r+\rho_1s}$ for the chosen index $i^*$ in Step (S3.3).  Hence, by Lemma~\ref{lemma-SCOBBA}(i),  $u_{i^*}^k- l_{i^*}^k> \frac{2\sqrt{\epsilon}}{\sqrt{r+\rho_1s}}$, i.e., the $i^*$-th edge of the sub-rectangle $\Delta^k$ must be longer than $\frac{2\sqrt{\epsilon}}{\sqrt{r+\rho_1s}}$.
According to Step (S3.3), it will be divided at either point $z^k_{i^*}$ or the midpoint $(u_{i^*}^k+l_{i^*}^k)/2$.
Note that if $u_{i}^k- l_{i}^k\le \frac{2\sqrt{\epsilon}}{\sqrt{r+\rho_1s}}$, then the $i$-th direction will never be chosen in Step (S3.3) as a branching direction at the $k$-iteration.
This means that all the edges of the sub-rectangle corresponding to a node with the smallest lower bound selected at every iteration will never be shorter than $\frac{2\sqrt{\epsilon}}{\sqrt{r+\rho_1s}}$.  Therefore, every edge
of the initial rectangle will be divided into at most $\left\lceil{\frac{\sqrt{r+\rho_1s}(z_u^i-z_l^i)}{2\sqrt{\epsilon}}}\right\rceil$ sub-intervals. In other words, to obtain an $\epsilon$-optimal solution  to Problem (\ref{OPDP-SD}), the total number of the relaxed subproblem (\ref{ODP:Bounded-CR}) required to be
solved  in all the runs of Algorithm~\ref{SCOBBA} is bounded above by
$$\prod_{i=1}^{r}\left\lceil{\frac{\sqrt{r+\rho_1s}(z_u^i-z_l^i)}{2\sqrt{\epsilon}}}\right\rceil.$$
This proves the conclusion of the theorem. \Halmos\endproof

\section{Tables in Section~\ref{sec:Experiment}}
\begin{table}[!ht]\caption{The statistics of the linear regression model for each stock in Example~\ref{exp_nasdaq}.}\label{tab:estimationstat}
\centering
\begin{tabular}{c|c|c|c|c|c|c}
\hline
stock symbol &GM &KO &PEP&WMT&AAPL& GE\\
\hline
p-value&0.000&0.000&0.000&0.000&0.000&0.000\\
\hline
R-square&51.10\% &63.89\% &56.41\% &40.09\% &47.78\% &48.09\%\\
\hline
\end{tabular}
\end{table}
 \begin{table}[!ht]\caption{Descriptive statistics of the twenty stocks in the sample month  considered in Example~\ref{exam2_nasdaq}.}\label{tab:example4-statistics}
 \centering
 \begin{tabular}{c|c|c|c}
 \hline
 \multirow{2}{*}{stock symbol} & aver. bid price & aver. hourly trade vol. &aver. trade size\\
 &(in dollar) &(in million) &(in share)\\
  \hline
 AMD&\$6.8335	&0.4651 &671.4899	
\\
 \hline
AXP&\$64.4741	&0.0544	&101.2464
\\
\hline
BAC&\$ 14.5604	&0.4974	&821.4500
\\
\hline
BLK&\$366.0428	&0.0126	 &61.3285

\\
\hline
BRK.B&\$144.2829	&0.0307	 &70.0429

\\
\hline
COST&\$167.1873	&0.0421	&64.8231

\\
\hline
 CVX&\$101.7120	&0.1135	&93.4919

 \\
 \hline
 DIS&\$96.0200	&0.1006	&100.6033

 \\
 \hline
 F&\$12.7467	&0.2529	&440.9061
\\
 \hline
 FB&\$124.9424	&0.6525 	&126.9192
\\
 \hline
 JNJ&\$125.0084	&0.0844	&86.7354
\\
 \hline
 MCD&\$118.0119	&0.1252	&94.7140
\\
 \hline
 MO&\$67.7486	&0.0791	&102.1315
\\
 \hline
 MSFT&\$56.5364	&0.6957	 &212.4948
\\
 \hline
 V&\$78.4346	&0.1048	&101.3619
\\
 \hline
 VZ&\$55.4707	&0.1976	&150.9227
\\
 \hline
 WFC&\$48.0434	&0.2071	&175.8796
\\
 \hline
 WMT&\$72.9371	&0.0704	&99.5293
\\
 \hline
 XOM&\$88.1901	&0.2000	&94.1738
\\
 \hline
YHOO&\$38.2431	&0.3934	&199.8870
\\
 \hline
 \end{tabular}
 \end{table}
 
 \section{The way to generate matrices with fixed number of negative eigenvalues in Example~\ref{exp_limitr}.}
To control the number of negative eigenvalues of the matrices $\hat{\Lambda}-\frac{1}{2}\hat{\Gamma}$ and $\hat{\Lambda}+\frac{1}{2}\hat{\Gamma}$, we generate the matrices $\Lambda$ and $\Gamma$ in the following way.
First, we randomly generate two matrices $C,F\in{\mathbb R}^{m\times m}$ whose entries are from $U[10^{-6},10^{-5}]$. Let $\tilde{C}=\frac{1}{2}(C+C^T)$ and $\tilde{F}=\frac{1}{2}(F+F^T)$. Then we compute the eigenvalues $\lambda_1\le\ldots\le\lambda_m$ and $\mu_1\le\ldots\le\mu_m$ for matrices $\tilde{C}$ and $\tilde{F}$, respectively. For arbitrary fixed integers $q$ and $s$ with $1\le q,~s<m$, let
\[v=\frac{1}{2}(\lambda_q+\lambda_{q+1}),~~w=\frac{1}{2}(\mu_s+\mu_{s+1}).\]
It is easy to check that matrices $\tilde{C}-vI$ and $\tilde{F}-wI$ have exactly $q$ and $s$ negative eigenvalues, respectively.
Let us define
  $$\Lambda=\frac{1}{2}(C-vI+F-wI),~~\Gamma=C-vI-F+wI. $$%~~\hat\Lambda=\frac{1}{2}(\Lambda+\Lambda^T),~~\hat\Gamma=\frac{1}{2}(\Gamma+\Gamma^T).$$
 It is then easy to check that  $\hat{\Lambda}-\frac{1}{2}\hat{\Gamma}=\tilde{F}-wI$ and $\hat{\Lambda}+\frac{1}{2}\hat{\Gamma}=\tilde{C}-vI$,
 where $\hat\Gamma$ and $\hat\Lambda$ are given in \eqref{eqn:sys-matrices}.

\end{document}